\documentclass[journal ]{new-aiaa}
\usepackage[utf8]{inputenc}
\usepackage{textcomp}
\usepackage{graphicx}
\usepackage{amsmath}
\usepackage[version=4]{mhchem}
\usepackage{siunitx}
\usepackage{longtable,tabularx}
\usepackage{algpseudocode}
\usepackage{algorithm}
\usepackage{multirow}
\usepackage{booktabs} 
\usepackage{caption}
\usepackage{subcaption}
\usepackage{makecell}
\usepackage{mathtools} 
\usepackage[dvipsnames]{xcolor}

\global\long\def\norm#1{\lVert#1\rVert}%
\global\long\def\abs#1{\lvert#1\rvert}%

\newcommand{\behcet}{Beh\c{c}et}
\newcommand{\acikmese}{A\c{c}\i{}kme\c{s}e}

\newcommand{\bartraj}{\{\bar{{\bf x}}_k, \bar{{\bf u}}_k,\bar{s}_k\}}
\newcommand{\hattraj}{\{\hat{{\bf x}}_k, \hat{{\bf u}}_k,\hat{s}_k\}}
\newcommand{\traj}{\{{{\bf x}}_k, {{\bf u}}_k,{s}_k\}}
\newcommand{\defeq}{\vcentcolon=}
\newcommand{\change}[1]{\textcolor{black}{#1}} 

\newtheorem{theorem}{Theorem}

\usepackage{tcolorbox}
\tcbuselibrary{skins}
\newtcolorbox{mybox}
{
  enhanced jigsaw,
  colframe=black,
  colback=white,
  drop shadow=black!50!white,
  boxrule=0.75pt
}

\newtcolorbox{mybox2}
{
  enhanced jigsaw,
  colframe=black,
  colback=white,
  drop shadow=black!50!white,
  boxrule=0.75pt,
  hbox
}

\title{Six-Degree-of-Freedom Aircraft Landing Trajectory Planning with Runway Alignment}

\author{Taewan Kim$^*$\!, Abhinav G. Kamath$^*$\!, Niyousha Rahimi\footnote{Ph.D.\ Student, William E.\ Boeing Department of Aeronautics \& Astronautics; \texttt{\{twankim,\,agkamath,\,nrahimi\}@uw.edu}}\!, \behcet{} \acikmese{}$^\dagger$\!, and Mehran Mesbahi\footnote{Professor, William E.\ Boeing Department of Aeronautics \& Astronautics; \texttt{\{behcet,\,mesbahi\}@uw.edu}}
}
\affil{University of Washington, Seattle, WA 98195, USA}
\author{Jasper Corleis\footnote{Senior Software Engineer; \texttt{Jasper.P.Corleis@boeing.com}}}
\affil{Boeing, Seattle, WA 98204, USA}

\begin{document}

\maketitle

\begin{abstract}
This paper presents a numerical optimization algorithm for generating approach and landing trajectories for a six-degree-of-freedom (6-DoF) aircraft. We improve on the existing research on aircraft landing trajectory generation by  formulating the trajectory optimization problem with additional  real-world  operational constraints, including 6-DoF aircraft dynamics, runway alignment, constant wind field, and obstacle avoidance, to obtain a continuous-time nonconvex optimal control problem. Particularly, the runway alignment constraint enforces the trajectory of the aircraft to be aligned with the runway only during the final approach phase. This is a novel feature that is essential for preventing an approach that is either too steep or too shallow. The proposed method models the runway alignment constraint through a multi-phase trajectory planning scheme, imposing alignment conditions exclusively during the final approach phase. We compare this formulation with the existing state-triggered constraint formulation for runway alignment. To solve the formulated problem, we design a novel sequential convex programming algorithm called xPTR that extends the penalized trust-region (PTR) algorithm by incorporating an extrapolation step to expedite convergence. We validate the proposed method through extensive numerical simulations, including a Monte Carlo study, to evaluate the robustness of the algorithm to varying initial conditions.
\end{abstract}

\section*{Nomenclature}

{\renewcommand\arraystretch{1.0}
\noindent\begin{longtable*}{@{}l @{\quad=\quad} l@{}}
\multicolumn{2}{@{}l}{Aircraft dynamics}\\
$m$ & aircraft mass \\
$\bf{v}$, $\dot{\bf v}$ & aircraft velocity vector and time-derivative of velocity vector \\
$\bf{\Omega}$, $\dot{\bf \Omega}$ & angular velocity vector and time-derivative of angular velocity vector \\
$F_e$, $F_a$, $F_g$ & engine force, aerodynamic force, and gravity force \\
$p$, $q$, $r$ & components of angular velocity along the $x$-, $y$-, and $z$-axes, respectively \\
$V$ & magnitude of aircraft velocity (speed) \\
$\phi$, $\theta$, $\psi$ & roll, pitch, and yaw angles \\
$R_\mathcal{I}^\mathcal{B}$, $R_\mathcal{B}^\mathcal{I}$ & rotation matrices between the inertial frame and the body-fixed frame \\
$p_N$, $p_E$, $p_D$, $\bf p$ & north, east, and down components of position, and the position vector \\
$J$ & moment of inertia matrix \\
$M_e$, $M_a$ & engine moment, aerodynamic moment \\
${\bf v}_a$ & airspeed vector \\
$V_a$ & magnitude of airspeed \\
$\alpha$, $\beta$ & angle of attack and sideslip angle \\
$C_L$, $C_D$, $C_Y$ & coefficients of lift, drag, and sideforce \\
$C_{L_{wb}}$, $C_{L_{t}}$ & lift coefficients acting on wing and body, and tail \\
$C_{L_\alpha}$, $C_{L_{\alpha_t}}$, $C_{L_{qV}}$ & stability derivatives for lift\\
$S$, $S_t$ & \change{platform area for the wing and the tail} \\
$\alpha_t$ & angle of attack at the tail \\
$\epsilon$, $\epsilon_\alpha$ & downwash angle and slope constant for downwash angle\\
$l_t$ & longitudinal distance between aerodynamic center of tailplane and that of wing and body \\
$C_{D0}$, $C_{D1}$, $C_{D2}$, $C_{D_\alpha}$  & stability derivatives for drag \\
$C_{Y\beta}$, $C_{Y \delta_R}$ & stability derivative and control derivative for sideforce \\
$\bar{c}$ & mean aerodynamic chord \\
$C_{l\beta}$, $C_{lp}$, $C_{lr}$ & stability derivatives for rolling moment \\
$C_{l\delta_{A}}$, $C_{l\delta_{R}}$ & control derivatives for rolling moment \\
$C_{m0}$, $C_{m\alpha}$, $C_{mq}$ & stability derivatives for pitching moment \\
$C_{m\delta_{E}}$ & control derivatives for pitching moment \\
$C_{n\beta}$, $C_{n\alpha\beta}$, $C_{np},C_{nr}$ & stability derivatives for yawing moment \\
$C_{n\delta_{A}}$, $C_{n\delta_{R}}$ & control derivatives for yawing moment \\
$L$, $D$, $Y$ & lift, drag, and sideforce \\
$\bar{q}$ & dynamic pressure \\
$\rho$ & air density \\
$\bar{L}$, $\bar{M}$, $\bar{N}$ & rolling, pitching, and yawing moments \\
${\bf p}_{cg}$, ${\bf p}_{ac}$, ${\bf p}_T$ & positions of center of gravity, aerodynamic center, and engine \\
$\delta_A$, $\delta_E$, $\delta_R$ & commands for aileron, elevator, rudder \\
$\delta_T$, $\eta_T$ & thrust command and actual thrust command \\
\multicolumn{2}{@{}l}{Subscripts}\\
$\mathcal{I}$ & the north-east-down frame (inertial frame) \\
$\mathcal{B}$ & the aircraft-body frame (body-fixed frame) \\
$k$	& discrete time index \\
\multicolumn{2}{@{}l}{Notation}\\
$\mathbb{R}$, $\mathbb{R}_{+}$, $\mathbb{R}_{++}$, $\text{\ensuremath{\mathbb{Z}}}$ & sets of real, nonnegative real, and positive real numbers, and integers \\
$\mathbb{R}^n$ & $n$-dimensional Euclidean space \\
$\mathcal{N}_q^r$ & a finite set of consecutive non-negative integers $\{q,q+1,\ldots,r\}$ \\
$\times$ & outer product operation \\
$\leq$ & component-wise inequality \\
$I_n$ & $n\times n$ identity matrix \\
$0_{n\times m}$ & $n\times m$ matrix of zeros \\
$\norm{\cdot}_2$ & vector 2 norm \\
$(a, b, c)$ & stacked vector of $a, b, c$, that is, $[a^\top, b^\top, c^\top]^\top$
\end{longtable*}}

\section{Introduction}
\lettrine{A}{pproach} and landing are crucial parts of flight since most aviation accidents happen during these phases, as indicated by a review of accident statistics \cite{airplanes2007statistical}.  Consequently, autonomous trajectory generation for aircraft approach and landing can contribute to flight safety and reliability. The major challenge in such trajectory generation schemes is that there are various operational constraints that must be taken into account, for instance the aircraft dynamics, the control limits, obstacle avoidance, and runway alignment. To tackle this, we employ an optimization-based trajectory generation method for the aircraft approach and landing that incorporates these operational constraints.

Optimization-based trajectory generation is a promising direction because it can provide a systematic framework to enforce state and control constraints, including the equations of motion (aircraft dynamics), i.e., it can generate dynamically feasible trajectories \cite{malyuta2022convex}. In addition to imposing constraints, another motivation for using an optimization-based approach is its capability in optimizing key mission objectives such as minimizing time and fuel. Due to these advantages, optimization-based approaches have been extensively utilized for various applications: powered descent guidance for reusable rockets \cite{acikmese2007convex,acikmese2013flight,reynolds2020dual,kamath2022real,kamath2023customized}, hypersonic reentry \cite{liu2016entry,mceowen2022hypersonic,mceowen2023high}, path generation and tracking for ground vehicles \cite{williams2016aggressive,kim2017path}, control of humanoid robots \cite{tassa2012synthesis,koenemann2015whole}, and spacecraft rendezvous \cite{malyuta2020fast,malyuta2023fast}.

Despite these advantages, it is challenging to design optimization-based methods to solve the six-degree-of-freedom (6-DoF) aircraft approach and landing trajectory generation. This is because, first, the trajectory generation problem needs to be expressed as a continuous-time (and hence, infinite-dimensional) optimal control problem over functions spaces, possibly having no analytical solutions. Second, the problem involves high-dimensional aircraft dynamics, considering nonlinear aerodynamic effects as constraints. Third, besides the aircraft dynamics, there are constraints that must be considered for real-world operational trajectory design. One such example is that of the runway alignment constraint that, as its name suggests, enforces the aircraft to be kept aligned with the runway centerline once the aircraft enters the final approach phase. This constraint is essential to generate trajectories for aircraft approach and landing in order to not only follow the standard landing procedure but also reduce any possibility of encountering any danger. Furthermore, we may need to impose obstacle avoidance constraints to avoid contrails or other aircraft.

In this paper, we formulate the aircraft approach and landing trajectory optimization problem to generate dynamically feasible trajectories, and provide a solution method based on sequential convex programming (SCP) that can reliably solve the formulated problem. For the problem formulation, we first model the 6-DoF aircraft dynamics using a research civil aircraft model (RCAM) that was developed by the Group for Aeronautical Research and Technology in Europe \cite{magni1997robust}, and then impose the modeled dynamics as a constraint to generate dynamically feasible trajectories. To impose the runway alignment constraint, we provide a multi-phase trajectory planning scheme \cite{kamath2022real} that separates the entire flight into two phases: the base leg and the final approach \cite{federal2011airplane}. Then, the alignment constraint is only imposed for the final approach phase, while preserving the convexity of the constraint and enabling exact temporal triggering. We compare the proposed multi-phase scheme-based approach with an approach that utilizes the form of state-triggered constraints (STCs) \cite{szmuk2020successive,reynolds2020dual} where the alignment condition is only enforced when the altitude of the aircraft is less than the user-specified value, resulting in a nonconvex constraint. For the solution method, we propose a novel SCP method, called the extrapolated penalized trust region method, or xPTR for short, that adopts and extends the penalized trust-region (PTR) method \cite{reynolds2020real} by employing an extrapolation update. Finally, we conduct extensive numerical simulations to validate the proposed method.

\subsection{Related works}\label{chap:related_work}
A number of previous studies have examined optimization-based aircraft trajectory generation \cite{kelley1971flight,burrows1983fuel,kato1986interpretation,miele1990optimal,lu1995optimal,delahaye2014mathematical}. The existing literature can be categorized according to the types of aircraft dynamics used. A two-dimensional horizontal planar model for example was used for commercial aircraft trajectory optimization by sequential quadratic programming (SQP), and dynamic programming in \cite{betts1995application} and \cite{hagelauer1998soft}. A similar planar model was used in \cite{sridhar2011aircraft}, where the authors constructed contrails formation models for obstacle avoidance trajectory planning. The work in \cite{vian1989trajectory} studied obstacle avoidance trajectory optimization with a point-mass 3-DoF aircraft model, in which the indirect approach based on Pontryagin's maximum principle (PMP) was used to solve the optimal control problem. On the other hand, the direct approach based on SQP was applied to generate 3-DoF aircraft trajectories \cite{raivio1996aircraft}. Trajectory generation for the 6-DoF aircraft model was also studied in \cite{desai2008six,bittner2012multi}. The work in \cite{desai2008six} employed a two-timescale collocation architecture to separate low-frequency translational motion and high-frequency rotational motion. In \cite{bittner2012multi}, to reduce the computation burden of trajectory optimization for the 6-DoF model, they first computed the optimal trajectory using the 3-DoF model, and then the obtained solution was used as an initial guess for the 6-DoF aircraft trajectory optimization. A flatness-based 6-DoF aircraft model was implemented for dynamic soaring trajectory optimization in \cite{elango2018trajectory,laad2020fourier}. This flatness-based approach was extended to incorporate the 6-DoF aircraft model in a model predictive control framework in \cite{sandeepkumar2022flatness}.

Aircraft approach and landing trajectory generation, in particular, has been also studied by researchers. In the trajectory generation for autoland systems, as described in \cite{siegel2011development,fallast2017automated}, Dubins paths are utilized, which employ relatively simple kinematic constraints compared to more complex dynamical models. The study in \cite{tsiotras2011initial} provided a framework for generating initial guesses for 3-DoF aircraft landing trajectory optimization. In \cite{zhao2010time,zhao2013time}, the authors proposed a method computing time-optimal state and input profiles, including a speed profile, a time-parameterized path, and a control history, to follow a prescribed geometric path for 3-DoF aircraft landing. In contrast, their study in \cite{zhao2013analysis} focused on generating energy-optimal state and input profiles for a given geometric path. Other researchers have also explored related topics. The authors in \cite{meng2014novel} proposed an algorithm to generate emergency landing trajectories for the 3-DoF aircraft model using a Gauss pseudospectral method, and \cite{hong2019adaptive} studied adaptive trajectory planning for impaired aircraft landing by integrating both parameter estimation and trajectory generation. In contrast to previous research, our study focuses on generating dynamical, optimal, and feasible approach and landing trajectories for aircraft, with a particular emphasis on incorporating 6-DoF aircraft dynamics under a constant wind field and operational constraints such as runway alignment and obstacle avoidance. With these additional considerations, we can generate more dynamic and reliable trajectories that are representative of real-world aircraft behavior. Unlike some previous studies that have focused on obtaining state and input profiles to follow prescribed geometric paths, our proposed method focuses on generating both state and input trajectories, including geometric paths.

\subsection{Contributions}

First, this paper provides a problem formulation for aircraft approach and landing trajectory generation. This formulation takes into account operational constraints such as 6-DoF aircraft dynamics, a runway alignment constraint, constant wind fields, and obstacle avoidance together, which have been left largely unaddressed in the literature. Particularly, we highlight the runway alignment constraint, ensuring that the resulting trajectory maintains an appropriate glide slope angle during the final approach phase.

Second, we propose an extension to the penalized trust region (PTR) algorithm for sequential convex programming \cite{reynolds2020real,szmuk2020successive}, called the extrapolated penalized trust region (xPTR) algorithm, to efficiently solve the formulated nonconvex aircraft approach and landing optimal control problem (OCP). The proposed algorithm has two distinct features that are absent in the published SCP-based application literature \cite{kamath2022real,malyuta2023fast,szmuk2020successive,reynolds2020dual,bonalli2019gusto}: (1) continuous-time constraint satisfaction and (2) the extrapolation update. To solve the continuous-time OCPs, numerical algorithms require them to be converted into finite-dimensional optimization problems by introducing discrete temporal node points. Consequently, constraint satisfaction is approximated by enforcing the constraints only at a finite number of temporal nodes \cite{szmuk2020successive}. Hence, this approximation could potentially cause inter-sample constraint violation \cite{dueri2017trajectory}. Our paper resolves this issue by integrating the integral formulation of the continuous-time constraints, as proposed in \cite{ctcs2024,teo1989nonlinear,lin2014control} and adopted in \cite{nmpc2024,chari2024fast}, into our SCP framework so that the method can prohibit inter-sample constraint violation between the discretization node points.

SCP algorithms are susceptible to the so-called \textit{crawling phenomenon} \cite{reynolds2020crawling}, where the progress of iterates slows down, possibly caused by the use of Lagrangian like function in the cost for the convex approximation of the nonconvex constraints and penalization of trust regions. To address this issue, we employ an extrapolation step, motivated by the Nesterov's acceleration technique \cite{nesterov1983method} and recently developed extrapolated first-order methods \cite{yang2017proximal,yu2022extrapolated}. Extrapolation aims to propel the solution further in the direction of the solution to the current convex subproblem with respect to the solution from the previous iteration—rather than merely accepting the solution to the subproblem at each iteration—which can accelerate convergence.

Finally, the outlined contributions of the proposed algorithm are validated via numerous numerical simulations. More specifically, comparisons are conducted between cases with and without the runway alignment constraint to highlight its necessity. Also, the proposed approach with the multi-phase scheme and the generic STC approach for imposing the runway alignment constraint are compared in terms of the quality of the generated trajectory and the convergence behavior. We also validate the proposed method in terms of continuous-time constraint satisfaction by means of a scenario wherein the resulting trajectory, without the proposed constraint reformulation, results in large constraint violations between temporal nodes, i.e., in the subintervals. A Monte Carlo study with various initial conditions is conducted to demonstrate the success-rate and the execution-time of the proposed algorithm. Within the study, we showcase convergence behavior in relation to the extrapolation parameter.

Based on these contributions, the proposed method can serve as a guidance module in autoland systems, providing reference trajectories to low-level tracking control modules. Given that the generated trajectory is dynamically feasible for the 6-DoF model, controllers will experience reduced error compared to trajectories computed using previous methods, as discussed in Section~\ref{chap:related_work}. Also, the proposed method can offer potential trajectory candidates for pilots to follow, accommodating varying flight conditions such as maximum speed, wind velocity, and locations of obstacles.

\subsection{Outline}

We present the problem statement in Section \ref{sec:problem_formul}, where the continuous-time optimal control problem is derived in Subsection \ref{subsec:cont_time_ocp}. Section \ref{sec:SCP} illustrates the details of the proposed xPTR algorithm, the subproblem template for which is provided in Subsection \ref{subsec:summary_discrete_OCP}. The numerical simulation results that validate the proposed method are provided in Section \ref{sec:simulation}. Finally, concluding remarks are provided in Section \ref{sec:conclusion}.



\section{Problem formulation} \label{sec:problem_formul}

In this section, we present a 6-DoF aircraft approach and landing trajectory generation problem with real-world operational constraints. 

\subsection{Assumptions} \label{subsec:assumption}

In this paper, we consider aircraft approach and landing maneuvers within a range of around 70 kilometers from the final landing site, with a corresponding flight horizon of around 10 minutes. Since the flight horizon and the distance are not too large, we assume that the effect of planetary rotation of the earth is negligible, and we use stationary atmosphere with constant air density, a flat non-rotating earth model, and a uniform gravitational field for the derivation of the aircraft dynamics. We assume a constant wind field without wind shear, and a rigid-body aircraft model, to simplify the system dynamics in the optimal control problem.

\subsection{6-DoF aircraft model} \label{subsec:aircraftmodel}

\begin{figure}
\centering
    \begin{mybox2}
        \centering
        \includegraphics[width=9cm]{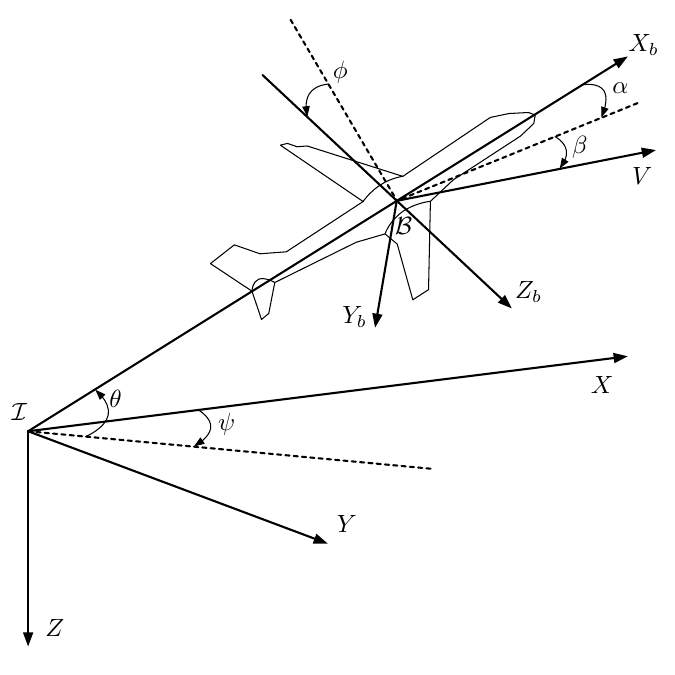}
    \end{mybox2}
\caption{Figure of 6-DoF aircraft coordinate frames}
\label{fig:aircraft_frame}
\end{figure}

We consider a 6-DoF aircraft model that is built upon a research civil aircraft model (RCAM) developed by the Group for Aeronautical Research and Technology in Europe \cite{magni1997robust}. Although the derivation of the aircraft dynamics itself is not a contribution of this paper, we include it here for the sake of completeness. More details in the derivation of the 6-DoF aircraft dynamics can be found in the following texts: \cite{stevens2015aircraft,etkin2012dynamics}. In the derivation, the time argument will be omitted if it is clear by context.

\subsubsection{Translational motion}
We define the inertial frame, denoted by subscript $\mathcal{I}$, to be in the North-East-Down (NED) coordinates and the body-fixed frame, denoted by subscript $\mathcal{B}$, to be in the aircraft body coordinates ($ABC$). These frames are illustrated in Figure \ref{fig:aircraft_frame}. The translational equations of motion, expressed in the body frame, can be derived from Newton's second law of motion as follows:
\begin{align}
F_{e}(t)+F_{a}(t)+F_{g}(t)=m(\dot{{\bf v}}(t)+{\bf \Omega}(t)\times{\bf v}(t)),\label{eq:dyn_trans}
\end{align}
where $m\in\mathbb{R}_{++}$ is the mass of the aircraft, ${\bf v}(t)\in\mathbb{R}^{3}$ is the airspeed of the aircraft, $\dot{{\bf v}}(t)\in\mathbb{R}^{3}$ is the time-derivative of the airspeed, ${\bf \Omega}(t)\in\mathbb{R}^{3}$ is the angular velocity, and $F_{e}(t),F_{a}(t),F_{g}(t)\in\mathbb{R}^{3}$ are the engine thrust, the aerodynamic force, and the gravitational force, respectively. These vectors in the body frame can be written component-wise as
${\bf v}=(u,v,w)$, $\dot{\bf v}=(\dot{u},\dot{v},\dot{w})$, and ${\bf \Omega}=(p,q,r)$. The magnitude of the airspeed is denoted as $V(t)\in\mathbb{R}$, where $V(t) \defeq \norm{{\bf v}(t)}_{2}$.

To express the airspeed vector ${\bf v}$ in the inertial frame, we need to perform a coordinate transformation from the body frame to the inertial frame. For this transformation, we use the roll $\phi(t)\in\mathbb{R}$, pitch $\theta(t)\in\mathbb{R}$, and yaw $\psi(t)\in\mathbb{R}$ Euler angles with the $ZYX$ convention. Then, the rotation matrix $R_{\mathcal{I}}^{\mathcal{B}}(\phi, \theta, \psi)\in\mathbb{R}^{3\times3}$ transforming the coordinates from the inertial frame to the body frame can be derived as follows
\begin{align*}
R(\phi, \theta, \psi)_{\mathcal{I}}^{\mathcal{B}} & =R_{\phi}R_{\theta}R_{\psi},\\
R_{\psi} & =\left[\begin{array}{ccc}
\cos\psi & \sin\psi & 0\\
-\sin\psi & \cos\psi & 0\\
0 & 0 & 1
\end{array}\right],\quad
R_{\theta}  =\left[\begin{array}{ccc}
\cos\theta & 0 & -\sin\theta\\
0 & 1 & 0\\
\sin\theta & 0 & \cos\theta
\end{array}\right],\quad
R_{\phi}  =\left[\begin{array}{ccc}
1 & 0 & 0\\
0 & \cos\phi & \sin\phi\\
0 & -\sin\phi & \cos\phi
\end{array}\right].
\end{align*}
It follows from $(R_{\mathcal{I}}^{\mathcal{B}})^{-1}=R_{\mathcal{B}}^{\mathcal{I}}$ that the time-derivative of position vector expressed in the inertia frame is
\begin{align}
\dot{{\bf p}}(t)=R_{\mathcal{B}}^{\mathcal{I}}(t){\bf v}(t) + {\bf w},\label{eq:dyn_NED}
\end{align}
where ${\bf p}(t)=(p_{N}(t),p_{E}(t),p_{D}(t))$ is a position vector consisting of the north, east, and down position components, respectively, and ${\bf w}(t)\in\mathbb{R}^3$ is the constant wind velocity with respect to the inertial frame. It is worth noting that the ground velocity of the vehicle with respect to the body frame is given by ${\bf v}_g = {\bf v} + R_{\mathcal{I}}^{\mathcal{B}}{\bf w}$.

\subsubsection{Rotational motion}

The rotational equations of motion expressed in the body frame can be derived from the Newton-Euler equations of motion as follows:
\begin{align}
M_{e}(t)+M_{a}(t)=J\dot{{\bf \Omega}}(t)+{\bf \Omega}(t)\times J{\bf \Omega}(t),\label{eq:dyn_rotation}
\end{align}
where ${\bf \Omega}(t)\in\mathbb{R}^{3}$ is the angular velocity vector that is expressed as ${\bf \Omega} \defeq (p, q, r)$ in the body frame.
The matrix $J\in\mathbb{R}^{3\times3}$ is the moment of inertia matrix that is expressed as
\begin{align*}
J=\left[\begin{array}{ccc}
J_{x} & 0 & -J_{xz}\\
0 & J_{y} & 0\\
-J_{xz} & 0 & J_{z}
\end{array}\right],
\end{align*}
where we assume that the aircraft is symmetric about the $xz$ plane. The moments $M_{e}(t),M_{a}(t)\in\mathbb{R}^{3}$ are generated by the engine thrust and the aerodynamic force, respectively. The relationship between Euler angles ${\bf \Phi}(t) \defeq (\phi,\theta,\psi)$ and the angular velocity ${\bf \Omega}(t)$ is given by
\begin{align}
\dot{{\bf \Phi}}(t) & =R^{E}(t){\bf \Omega}(t),\label{eq:dyn_euler}\\
R^{E} & \defeq\left[\begin{array}{ccc}
1 & \sin\phi\tan\theta & \cos\phi\tan\theta\\
0 & \cos\phi & -\sin\phi\\
0 & \sin\phi\sec\theta & \cos\phi\sec\theta
\end{array}\right].\nonumber 
\end{align}


\subsubsection{Aerodynamic forces and moments}

The aerodynamic forces and moments depend on the angle of attack $\alpha(t)\in\mathbb{R}$ and the sideslip angle $\beta(t)\in\mathbb{R}$ that have the following relationships with the airspeed:
\begin{align}
\tan\alpha & =\frac{w}{u},\quad\sin\beta=\frac{v}{V}.\label{eq:def_alpha_beta}
\end{align}
The angle of attack and sideslip angle are illustrated in Figure \ref{fig:aircraft_frame}.

The aerodynamic forces and moments are computed based on the aerodynamic coefficients for drag $C_{D}(t)\in\mathbb{R}$, sideforce $C_{Y}(t)\in\mathbb{R}$, and lift $C_{L}(t)\in\mathbb{R}$, that are functions of the angle of attack $\alpha$, the sideslip angle $\beta$, angular velocity ${\bf \Omega}$, and the control surfaces. The lift coefficient $C_{L}$ is given as
\begin{align*}
C_{L}=C_{L_{wb}}+C_{L_{t}},
\end{align*}
where $C_{L_{wb}}(t)\in\mathbb{R}$ is the lift coefficient at the wing and body and $C_{L_{t}}(t)\in\mathbb{R}$ is the lift coefficient at the tail. The coefficient $C_{L_{wb}}$ is assumed to have the following linear relationship with the angle of attack:
\begin{align*}
C_{L_{wb}}=C_{L_{\alpha}}(\alpha-\alpha_{0}),
\end{align*}
where $C_{L_{\alpha}}\in\mathbb{R}$ is the stability derivative parameter for $\alpha$, and $\alpha_{0}\in\mathbb{R}$ is the constant value of the angle of attack at which the lift for the wing and body becomes zero. The lift coefficient at the tail, $C_{L_{t}}$, is given as
\begin{align*}
C_{L_{t}}=C_{L_{\alpha_{t}}}\frac{S_{t}}{S}\text{\ensuremath{\alpha_{t}}},
\end{align*}
where $C_{L_{\alpha_{t}}}\in\mathbb{R}$ is the constant parameter, and $S_{t}\in\mathbb{R}$ and $S\in\mathbb{R}$ are the platform area for the tail and the wing, respectively. The angle of attack of the tail $\alpha_{t}(t)\in\mathbb{R}$ can be given by
\begin{align*}
\alpha_{t} & =\alpha-\epsilon+\delta_{E}+C_{L_{qV}}\frac{ql_{t}}{V},\\
\epsilon & =\epsilon_{\alpha}(\alpha-\alpha_{0}),
\end{align*}
where $\epsilon(t)\in\mathbb{R}$ is the downwash angle, $\epsilon_{\alpha}\in\mathbb{R}$ is the slope constant for the downwash angle, $\delta_{E}(t)\in\mathbb{R}$ is the elevator command, $C_{L_{qV}}\in\mathbb{R}$ is the stability derivative for the $y$-axis angular velocity $q$ and the magnitude of airspeed $V$, and $l_{t}\in\mathbb{R}$ is the longitudinal distance between the aerodynamic center of the tailplane and that of the wing and body.

The drag coefficient $C_{D}(t)\in\mathbb{R}$, which is a function of the angle of attack $\alpha$, is given as
\begin{align*}
C_{D}=C_{D0}+C_{D2}(C_{D_{\alpha}}\alpha+C_{D1})^{2},
\end{align*}
where $C_{D0}$, $C_{D1}$, $C_{D2}$, $C_{D_{\alpha}}\in\mathbb{R}$ are the stability derivative parameters related to how the drag changes. The sidefore coefficient $C_{Y}(t)\in\mathbb{R}$ is given as
\begin{align*}
C_{Y}=C_{Y\beta}\beta+C_{Y\delta_{R}}\delta_{R},
\end{align*}
where $C_{Y\beta}$, $C_{Y\delta_{R}}\in\mathbb{R}$ are the stability derivative and the control derivative parameters, respectively.

The aerodynamic moment coefficients for the lift are given by
\begin{align*}
C_{l} & =C_{l\beta}\beta+C_{lp}\frac{p\bar{c}}{V}+C_{lr}\frac{r\bar{c}}{V}+C_{l\delta_{A}}\delta_{A}+C_{l\delta_{R}}\delta_{R},\\
C_{m} & =C_{m0}+C_{m\alpha}\frac{S_{t}l_{t}}{S\bar{c}}(\alpha-\epsilon)+C_{mq}\frac{S_{t}l_{t}^{2}}{S\bar{c}}\frac{q}{V}+C_{m\delta_{E}}\frac{S_{t}l_{t}}{S\bar{c}}\delta_{E},\\
C_{n} & =C_{n\beta}\beta+C_{n\alpha\beta}\alpha\beta+C_{np}\frac{\bar{c}p}{V}+C_{nr}\frac{\bar{c}r}{V}+C_{n\delta_{R}}\delta_{R},
\end{align*}
where $\bar{c}\in\mathbb{R}_{++}$ is mean aerodynamic chord. The stability derivatives $C_{l\beta},C_{lp},C_{lr}\in\mathbb{R}$ and control derivatives $C_{l\delta_{A}},C_{l\delta_{R}}\in\mathbb{R}$ are for the rolling moment, and $C_{m0},C_{m\alpha},C_{mq}\in\mathbb{R}$ and $C_{m\delta_{E}}\in\mathbb{R}$ are the stability derivatives and the control derivative, respectively, for the pitching moment. For the yawing moment, the stability derivatives $C_{n\beta},C_{n\alpha\beta},C_{np},C_{nr}\in\mathbb{R}$ and the control derivatives $C_{n\delta_{A}},C_{n\delta_{R}}\in\mathbb{R}$ are used. More details regarding the effect of these stability and control derivatives on the aerodynamic moments can be found in \cite{magni1997robust,stevens2015aircraft}.

Given the coefficients $C_D$, $C_L$, and $C_Y$, the drag force $D(t)\in\mathbb{R}$, the lift force $L(t)\in\mathbb{R}$, and the sideslip force $Y(t)\in\mathbb{R}$ are derived as
\begin{align*}
D & =\bar{q}SC_{D},\\
L & =\bar{q}SC_{L},\\
Y & =\bar{q}SC_{Y},
\end{align*}
where $\bar{q}(t)\in\mathbb{R}_{+}$ is the dynamic pressure, which is given by
\begin{align*}
\bar{q}=\frac{1}{2}\rho V_{A}^{2}.
\end{align*}
The parameter $\rho\in\mathbb{R}_{++}$ is air density, which is assumed to be constant. The aerodynamic forces $F_{a}(t)\in\mathbb{R}^{3}$ with respect to the body frame can be expressed as
\begin{align*}
F_{a} & =\left[\begin{array}{c}
L\sin\alpha-D\cos\alpha\cos\beta-Y\cos\alpha\sin\beta\\
-D\sin\beta+Y\cos\beta\\
-L\cos\alpha-D\sin\alpha\cos\beta-Y\sin\alpha\sin\beta
\end{array}\right].
\end{align*}
The aerodynamic moments about the aerodynamic center are given by
\begin{align*}
\bar{L} & =\bar{q}S\bar{c}C_{l},\\
\bar{M} & =\bar{q}S\bar{c}C_{m}\\
\bar{N} & =\bar{q}S\bar{c}C_{n},
\end{align*}
where $\bar{L}(t),\bar{M}(t),\bar{N}(t)\in\mathbb{R}^{3}$ are the rolling, pitching, and yawing moments, respectively. The net aerodynamic moment $M_a\in\mathbb{R}^3$ acting on the aircraft at the center of gravity in \eqref{eq:dyn_rotation} can be given by
\begin{align*}
M_{a}=( \bar{L},\bar{M},\bar{N}) + F_{a}\times({\bf p}_{cg}-{\bf p}_{ac}),
\end{align*}
where ${\bf p}_{cg}\in\mathbb{R}^{3}$ is the position of center of gravity and ${\bf p}_{ac}\in\mathbb{R}^{3}$ is the position of aerodynamic center.

\subsubsection{Gravity, engine forces and moments}

The gravitational force acting on the aircraft \eqref{eq:dyn_trans} with respect to the body frame is given as
\begin{align*}
F_{g}=R_{\mathcal{I}}^{\mathcal{B}}(0,0,mg),
\end{align*}
where $g\in\mathbb{R}_{++}$ is the acceleration due to gravity. Note that the gravitational force does not generate any moment since it is acting at the center of the gravity of the aircraft. The net force generated by the engine thrust is
\begin{align*}
F_{e} = (2\delta_{T}mg,0,0)
\end{align*}
where $\delta_{T}(t)\in\mathbb{R}_{+}$ is the throttle. The net moment generated by the engine thrust is
\begin{align*}
M_{e}=({\bf p}_{cg}-{\bf p}_{T})\times F_{e},
\end{align*}
where ${\bf p}_{T}\in\mathbb{R}^{3}$ is the position of the engine.

\subsubsection{Engine model}

The actual aircraft typically has a delay between the thrust command provided and the actual thrust produced. To account for this behavior of the engine, we model the engine as a first order system as follows:
\begin{align}
\dot{\delta}_{T}=\frac{1}{\tau_{T}}(\eta_{T}-\delta_{T}),\label{eq:dyn_engine_model}
\end{align}
where $\eta_{T}(t)\in\mathbb{R}$ is the thrust command and $\tau_{T}\in\mathbb{R}_{++}$ is a time constant.

\subsubsection{State-space form}

By merging the equations of motion \eqref{eq:dyn_trans},\eqref{eq:dyn_NED},\eqref{eq:dyn_rotation},\eqref{eq:dyn_euler}, and \eqref{eq:dyn_engine_model}, we can express the 6-DoF aircraft model in the state-space form as follows:
\begin{align}
\dot{{\bf x}}(t)={\bf f}({\bf x}(t),{\bf u}(t)),\label{eq:dynamics_statespace}
\end{align}
where the state vector ${\bf x}(\cdot)\in\mathbb{R}^{13}$ and the input vector ${\bf u}(\cdot)\in\mathbb{R}^{4}$ are
\begin{align*}
{\bf x}&=[{\bf p}^{\top},{\bf v}^{\top},{\bf \Phi}^{\top},{\bf \Omega}^{\top},\delta_{T}]^\top, \\
{\bf u}&=[\delta_{A},\delta_{E},\delta_{R},\eta_{T}]^\top.
\end{align*}

\subsection{State and input constraints}

In this subsection, we describe the state and input constraints. First, the altitude of the aircraft must always be positive to prohibit the aircraft from colliding with the ground—this constraint is written as
\begin{align}
-p_{D}\geq0.\label{eq:const_height}
\end{align}
The minimum and maximum velocity constraints are given by
\begin{align}
{\bf v}_{\min}\leq{\bf v}\leq{\bf v}_{\max},\label{eq:const_velocity}
\end{align}
where the vectors ${\bf v}_{\min},{\bf v}_{\max}\in\mathbb{R}^{3}$ are the component-wise minimum and maximum velocity. The Euler angles and the angular velocities have the following constraints:
\begin{align}
\phi_{\min} \leq\phi\leq\phi_{\max},\quad\theta_{\min}\leq\theta\leq\theta_{\max},\quad{\bf \Omega}_{\min} \leq{\bf \Omega}\leq{\bf \Omega}_{\max},\label{eq:const_attitude}
\end{align}
where $\phi_{\min}$, $\theta_{\min}\in\mathbb{R}$, and ${\bf \Omega}_{\min}\in\mathbb{R}^{3}$ are the minimum allowable roll angle, pitch angle, and angular velocity, respectively, and $\phi_{\max}$, $\theta_{\max}\in\mathbb{R}$, and ${\bf \Omega}_{\max}\in\mathbb{R}^{3}$ are the maximum allowable roll angle, pitch angle, and angular velocity, respectively. To prevent the wings from stalling and having a negative lift force, we impose the upper and lower bounds on the angle of attack as follows:
\begin{align*}
\alpha_{\min}\leq\alpha\leq\alpha_{\max},
\end{align*}
where $\alpha_{\min},\alpha_{\max}\in[-\frac{\pi}{2},\frac{\pi}{2}]$. These constraints are equivalent to 
\begin{align}
u_{a}\tan\alpha_{\min}\leq w_{a} & \leq u_{a}\tan\alpha_{\max},\label{eq:const_alpha}
\end{align}
by the definition of $\alpha$ given in \eqref{eq:def_alpha_beta}. We impose the following box constraint on the input ${\bf u}$:
\begin{align}
{\bf u}_{\min}\leq{\bf u} & \leq{\bf u}_{\max},\label{eq:const_input}
\end{align}
where ${\bf u}_{\min},{\bf u}_{\max}\in\mathbb{R}^{4}$ are the minimum and maximum allowable control inputs, respectively. We impose the following constraints on the thrust:
\begin{subequations}
\label{eq:const_input_dot}
\begin{align}
\delta_{T,\min}&\leq\delta_{T}\leq\delta_{T,\max}, \\
\dot{\delta}_{\min}&\leq\frac{1}{\tau_{T}}(\eta_{T}-\delta_{T})\leq\dot{\delta}_{\max}, 
\end{align}
\end{subequations}
where $\delta_{T,\min}\in\mathbb{R}$ and $\delta_{T,\max}\in\mathbb{R}$ are the minimum and maximum throttle, respectively, and $\dot{\delta}_{\min}\in\mathbb{R}$ and $\dot{\delta}_{\max}\in\mathbb{R}$ are the minimum and maximum throttle rates, respectively.

\subsection{Boundary conditions}
The initial conditions are given by
\begin{subequations}\label{eq:boundary_initial}
\begin{align}
{\bf p}(t_{0}) & ={\bf p}_{i},\quad{\bf v}(t_{0})={\bf v}_{i}, \\
{\bf \Phi}(t_{0}) & ={\bf \Phi}_{i},\quad{\bf \Omega}(t_{0})={\bf \Omega}_{i},
\end{align}
\end{subequations}
where $t_{0}\in\mathbb{R}_+$ is the start time, and ${\bf p}_{i}\in\mathbb{R}^{3}$, ${\bf v}_{i}\in\mathbb{R}^{3}$, ${\bf \Phi}_{i}\in\mathbb{R}^{3}$, and ${\bf \Omega}_{i}\in\mathbb{R}^{3}$ are the initial position, velocity, Euler angles, and angular velocity, respectively. The terminal boundary conditions are given by
\begin{subequations}\label{eq:boundary_final}
\begin{align}
{\bf p}(t_{f}) & =0_{3\times1},\quad \abs{{\bf v}(t_{f})}\leq{\bf v}_{f}, \\
{\bf \Phi}(t_{f}) & =0_{3\times1},\quad{\bf \Omega}(t_{f})=0_{3\times1},
\end{align}
\end{subequations}
where $t_{f}\in\mathbb{R}_{++}$ is the final time. The final position at which the aircraft touches down is set to the origin of inertial frame $\mathcal{I}$. The airspeed ${\bf v}$ at the final time should be less than ${\bf v}_{f}\in\mathbb{R}_{+}^{3}$ component-wise to make sure that the velocity of the aircraft on the runway is not too large. 

\subsection{Runway alignment constraint}

As the aircraft descends towards the touch down point, it needs to be aligned with the runway in order to prevent an approach that is either too steep or too shallow for passenger comfort. This alignment constraint can be described using simple linear constraints that ensure that the position of the aircraft is aligned with the runway during the final approach phase. An important consideration here is that the switching time, which determines when the alignment constraint is activated and when the final approach phase starts, needs to be optimized alongside the entire trajectory. Hence, this constraint differs from the aforementioned pointwise state and input constraints that are enforced at all times during the entire flight.

\begin{figure}
\centering
    \begin{mybox2}
        \centering
        \includegraphics[scale=0.3]{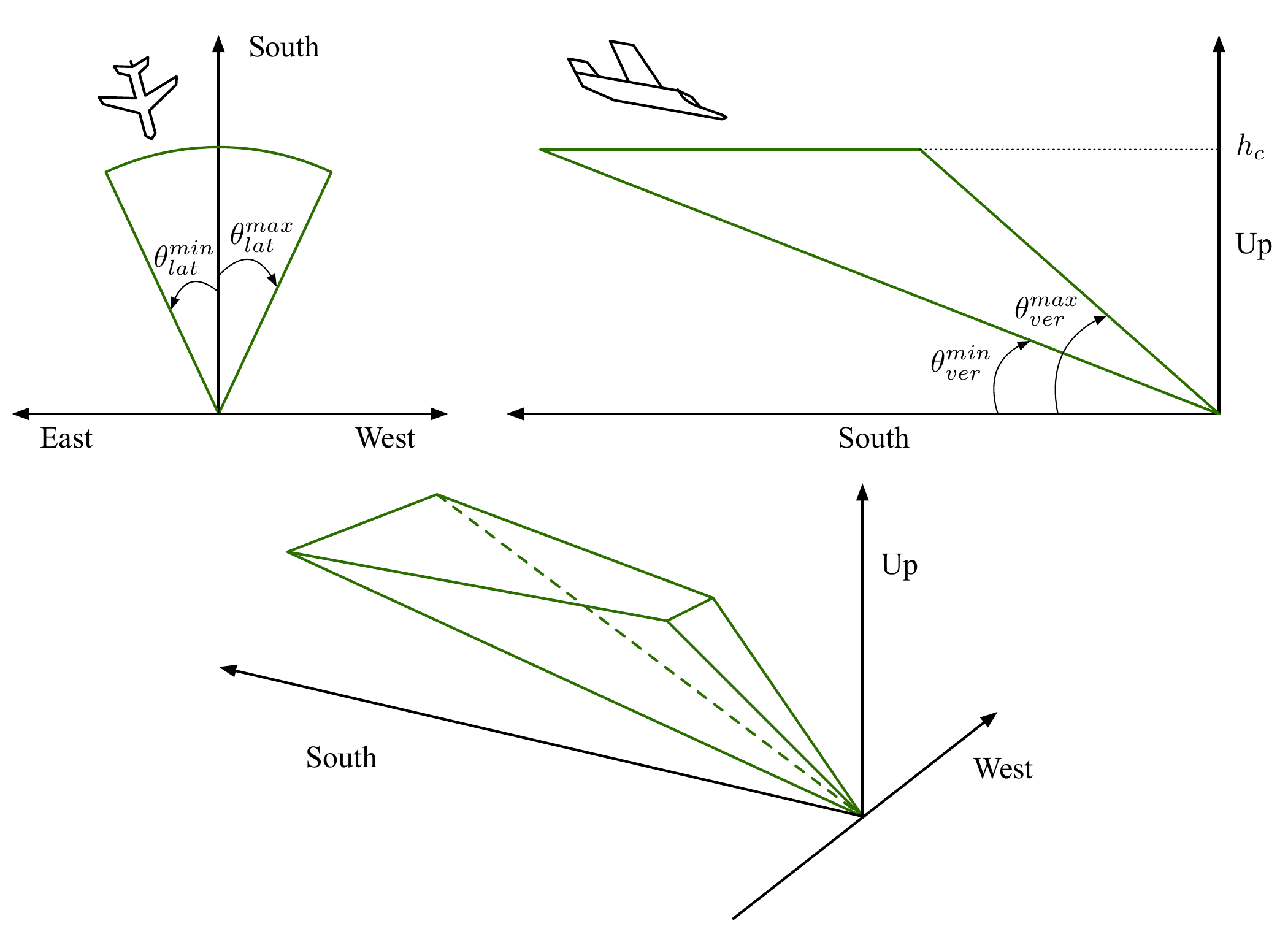}
    \end{mybox2}
\caption{Illustration of the runway alignment constraint.}
\label{fig:runway_alignment}
\end{figure}

The runway alignment constraint with the switching time $t_s$ can be mathematically stated as
\begin{align}
    \exists t_s\in [t_0, t_f] \text{ such that }-p_{D}(t)\leq h_{c} \text{ and }{\bf c}({\bf p}(t))\leq0, \quad \forall t\in [t_s, t_f],\label{eq:runway_align_time}
\end{align}
where $h_{c}\in\mathbb{R}_{++}$ is the user-specified altitude. The alignment constraint function ${\bf c}:\mathbb{R}^{3}\rightarrow\mathbb{R}^{4}$ is defined as
\begin{align}
\left[\begin{array}{c}
-p_{E}+p_{N}\tan\theta_{\mathrm{lat}}^{\max}\\
p_{E}-p_{N}\tan\theta_{\mathrm{lat}}^{\min}\\
-p_{D}+p_{N}\tan\theta_{\mathrm{ver}}^{\max}\\
p_{D}-p_{N}\tan\theta_{\mathrm{ver}}^{\min}
\end{array}\right]\leq0,\label{eq:constraint_in_align}
\end{align}
where angles $\theta_{\mathrm{lat}},\theta_{\mathrm{ver}}\in[-\frac{\pi}{2},\frac{\pi}{2}]$ are
\begin{align}
\theta_{\mathrm{lat}}=\tan^{-1}\left(\frac{-p_{E}}{-p_{N}}\right),\quad\theta_{\mathrm{ver}}=\tan^{-1}\left(\frac{-p_{D}}{-p_{N}}\right),\label{eq:def_angles}
\end{align}
representing the angles between the runway and the projected positions of the aircraft onto the $xy$-plane and the $xz$-plane, respectively, in the inertial frame. The parameters $\theta_{\mathrm{lat}}^{\min},\theta_{\mathrm{ver}}^{\min}\in\mathbb{R}$ and $\theta_{\mathrm{lat}}^{\max},\theta_{\mathrm{ver}}^{\max}\in\mathbb{R}$ are the minimum and maximum values of $\theta_{\mathrm{lat}}$ and $\theta_{\mathrm{ver}}$, respectively. The feasible area for the runway alignment constraint is illustrated in Figure \ref{fig:runway_alignment}.

The alignment constraint in \eqref{eq:runway_align_time} showcases an instance of existential quantification in the form of ``there exists a time $t_s$''. The existential quantifier can complicate the optimization problem since it does not specify a precise switching time $t_s$, but instead it indicates that such a time $t_s$ exists within the given interval \([t_0, t_f]\). In the next section, we discuss how this constraint can be incorporated into our SCP framework by using time-interval dilation variables.

An alternative way to express the runway alignment constraint is to utilize the state-triggered constraint (STC) formulation \cite{szmuk2020successive,reynolds2020dual}. The alignment constraint can then be given as the following logical statement:
\begin{align}
\text{if }-p_{D} < h_{c}\text{, then }{\bf c}({\bf p})\leq0.\label{eq:runway_align_STC}
\end{align}
The statement \eqref{eq:runway_align_STC} means that if the aircraft satisfies a trigger condition that is $-p_{D} < h_{c}$, the aircraft should align with the runway by satisfying the condition ${\bf c}({\bf p})\leq0$. On the other hand, if the trigger condition is not activated, which means the altitude is greater than $h_c$, the alignment condition is not necessarily satisfied. The alignment constraint given by STC \eqref{eq:runway_align_STC} is not equivalent to \eqref{eq:runway_align_time}. As an example, having both $-p_D < h_c$ and ${\bf c}({\bf p})>0$ before the switch time $t_s$ does not violate the original condition \eqref{eq:runway_align_time}, whereas the same condition violates \eqref{eq:runway_align_STC}. Hence, the alignment constraint in the STC form is more conservative. However, if we assume that the height of the altitude is monotonically decreasing, the two conditions \eqref{eq:runway_align_time} and \eqref{eq:runway_align_STC} are equivalent. 

The expression for the STC given in \eqref{eq:runway_align_STC} is a logical statement, and consequently cannot be readily imposed in an optimization setting. One way to impose the STC \eqref{eq:runway_align_STC} is to use a $\min$ function \cite{szmuk2020successive,reynolds2020dual,kamath2023customized} as follows:
\begin{align}
{\bf g({\bf p})} \coloneqq -\min(-p_D - h_c, 0)\,{\bf c}({\bf p}) \leq 0,\label{eq:runway_align_minfunction}
\end{align}
where the function ${\bf g}:\mathbb{R}^3 \rightarrow \mathbb{R}^4$ is introduced to compactly denote the constraint. One can find that the condition \eqref{eq:runway_align_minfunction} is equivalent to the original STC condition \eqref{eq:runway_align_STC}. One drawback of this expression is that \eqref{eq:runway_align_minfunction}  is nonconvex, which can in turn increase the problem complexity. 

\subsection{Obstacle avoidance constraint}

We consider obstacle avoidance by imposing additional constraints. The obstacles could be other aircraft or areas where the weather is so severe that it could be detrimental to flight safety, for example. The obstacles are modeled as three dimensional ellipsoids, and the obstacle avoidance constraints can be written as
\begin{align}
1-\norm{H_{i}({\bf p}(t)-r_{j})}_{2}\leq0,\quad\forall i\in\mathcal{N}_{1}^{n_{\mathrm{obs}}},\quad\forall t\in[t_{0},t_{f}],\label{eq:const_obstacle_avoidance}
\end{align}
where $H_{i}\in\mathbb{R}^{3\times3}$ is the shape matrix for the ellipse, $r_{i}\in\mathbb{R}^{3}$ is the center of the ellipse, $n_{\mathrm{obs}}\in\mathbb{Z}$ is the number of obstacles, where subscript $i$ represents the index of the obstacle under consideration.

\subsection{Performance Metrics}

The aircraft landing trajectory has two main performance indices: minimum-time and minimum-thrust. The cost function for each index can be formally written as
\begin{align}
J_{\text{time}}=t_{f},\quad J_{\mathrm{thr}}=\int_{t_{0}}^{t_{f}}\eta_{T}^{2}(t)\text{d}t,\label{eq:cost_time_thrust_cont}
\end{align}
respectively. We levy a penalty on the rate of the control input commands $\delta_{A}$, $\delta_{E}$, and $\delta_{R}$ to ensure that the resulting control input profiles are smooth. The cost function for this penalty is given as 
\begin{align}
J_{r}=\int_{t_{0}}^{t_{f}}\dot{\delta}_{A}(t)^{2}+\dot{\delta}_{E}(t)^{2}+\dot{\delta}_{R}(t)^{2}\text{d}t.\label{eq:cost_time_control_rate_cont}
\end{align}
Also, we penalize the angular velocity to prohibit rapid changes in the attitude of the aircraft. The corresponding cost function is given as
\begin{align}
J_{\Omega}=\int_{t_{0}}^{t_{f}}p(t)^{2}+q(t)^{2}+r(t)^{2}\text{d}t.\label{eq:cost_angular_vel_cont}
\end{align}
Hence, the combined objective function for the aircraft landing problem obtained by merging \eqref{eq:cost_time_thrust_cont}-\eqref{eq:cost_angular_vel_cont} can be given by 
\begin{align}
J=w_{t}J_{\text{time}}+w_{\mathrm{thr}}J_{\mathrm{thr}}+w_{r}J_{r}+w_{\Omega}J_{\Omega},\label{eq:cost_all_cont}
\end{align}
where $w_{t},w_{\mathrm{thr}},w_{r},w_{\Omega}\in\mathbb{R}$ are the weighting factors for the final time, the thrust, the rate of the control input commands, and the angular velocity, respectively. We are now ready to derive the problem formulation provided in Section \ref{subsec:cont_time_ocp}.

\begin{mybox}
    \vspace{-0.75em}
    \subsection{The continuous-time nonconvex aircraft landing trajectory optimization problem} \label{subsec:cont_time_ocp} 
    \vspace{-2.5em}
    \begin{align*}
        \underset{t_{f},{\bf x}(\cdot),{\bf u}(\cdot)}{\text{minimize}}
         \quad& \text{Eq \eqref{eq:cost_all_cont}} \\
        \text{subject to} \quad& \forall t\in[t_{0},t_{f}],\nonumber \\
        \fbox{\text{6-DoF dynamics}} \quad& \text{Eq \eqref{eq:dynamics_statespace}} \\
        \fbox{\text{State constraints}} \quad& \text{Eqs \eqref{eq:const_height}-\eqref{eq:const_alpha}} \\
         \fbox{\text{Input constraints}} \quad& \text{Eqs \eqref{eq:const_input}-\eqref{eq:const_input_dot}} \\
        \fbox{\text{Runway alignment}} \quad& \text{Eq \eqref{eq:runway_align_time} or  \eqref{eq:runway_align_minfunction}} \\
        \fbox{\text{Obstacle avoidance}} \quad& \text{Eq \eqref{eq:const_obstacle_avoidance}} \\
        \fbox{\text{Boundary conditions}} \quad& \text{Eqs \eqref{eq:boundary_initial}-\eqref{eq:boundary_final}}
    \end{align*}
    \vspace{-2.5em}
\end{mybox}

\section{Sequential Convex Programming} \label{sec:SCP}

We design an SCP-based algorithm to solve the free-final-time continuous-time nonconvex optimal control problem derived in Section \ref{subsec:cont_time_ocp}. The core idea in SCP is to sequentially construct fixed-final-time discrete-time convex subproblems that approximate the original problem. The algorithm involves iteratively solving these subproblems until it converges to a solution to Section \ref{subsec:cont_time_ocp}, wherein the nonconvex constraints at the current iteration are linearized about the solution from the previous iteration. In this paper, we develop the xPTR algorithm by combining the PTR algorithm with an extrapolation update. The xPTR algorithm employs multiple shooting discretization, in which the continuous-time input signal is represented by a finite number of parameters. Further, a time scaling transformation known as time-interval dilation is adopted. With discretization and dilation, it becomes possible to convert the free-final-time continuous-time problem into a fixed-final-time discrete-time problem. Then, to construct the convex subproblems approximating the original nonconvex problem, all nonconvex constraints including the system dynamics are linearized around the solution to the previous subproblem, i.e., around the reference trajectory. 
Given the solution to this convex subproblem, the algorithm performs an extrapolation to update the reference trajectory. The entire procedure is repeated until convergence, at which the solution is dynamically feasible, i.e., the nonlinear dynamics are exactly satisfied up to integration tolerances.

\subsection{Multiple shooting discretization and time-scaling transformation}

\begin{figure}
    \centering

    \begin{mybox2}
        \centering
        \includegraphics[scale=0.6]{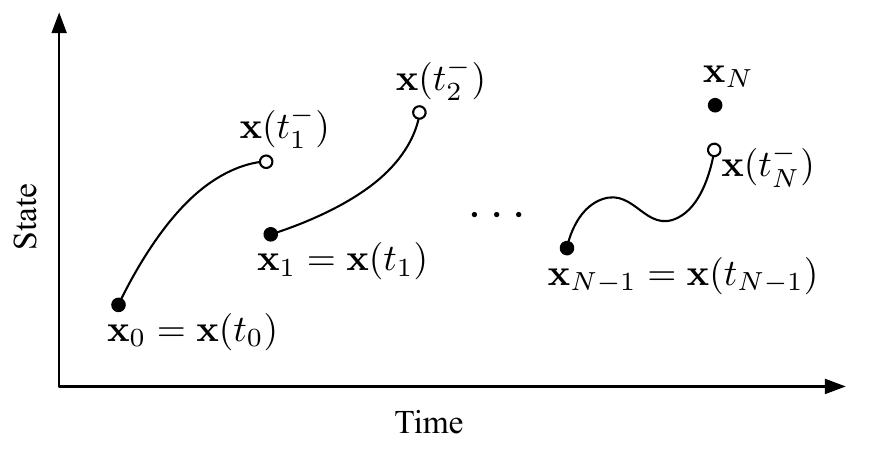}
    \end{mybox2}
\caption{Illustration of multiple shooting discretization.}
\label{fig:multiple_shooting}
\end{figure}

To numerically solve the given continuous-time optimal control problem in Section \ref{subsec:cont_time_ocp}, we first perform multiple shooting discretization, which was first proposed in \cite{bock1984multiple}. 
We parameterize the control input, choose a time mesh, and define the state and input nodes as follows:
\begin{align*}
t_0 &< t_1 < \cdots < t_k < \cdots < t_N = t_f,\quad\forall k\in\mathcal{N}_{0}^{N},\\
{\bf x}_k &\coloneqq {\bf x}(t_k), \quad {\bf u}_k \coloneqq {\bf u}(t_k),
\end{align*}
where $N\in\mathbb{R}_+$ is the number of subintervals.
Then, the control input signal ${\bf u}(\cdot)$ is parameterized using a first-order-hold (FOH) (piecewise continuous and affine) interpolation between successive input nodes as follows:
\begin{subequations}
\label{eq:input_FOH}
\begin{align}
{\bf {\bf u}}(t) & =\lambda_{k}^{m}(t){\bf u}_{k}+\lambda_{k}^{p}(t){\bf u}_{k+1},\quad\forall t\in[t_{k},t_{k+1}],\\
\lambda_{k}^{m}(t) & =\frac{t_{k+1}-t}{t_{k+1}-t_{k}},\quad\lambda_{k}^{p}(t)=\frac{t-t_{k}}{t_{k+1}-t_{k}}. 
\end{align}
\end{subequations}
Then, the state signal ${\bf x}(\cdot)$ in each subinterval is given by the solution to the following integral equation:
\begin{align}
{\bf x}(t) & ={\bf x}_{k}+\int_{t_{k}}^{t}{\bf f}({\bf x}(s),{\bf u}(s))\text{d}s,\quad \forall t\in[t_{k},t_{k+1}). \label{eq:integration_equation}
\end{align}
Figure \ref{fig:multiple_shooting} illustrates the relationship between ${\bf x}_{k}$ and ${\bf x}(t)$. The continuous-time state trajectory ${\bf x}(t)$ for each subinterval $t\in[t_{k},t_{k+1})$ is computed by solving \eqref{eq:integration_equation} from the point ${\bf x}_{k}$ at $t=t_{k}$ with the input signal ${\bf {\bf u}}(t)$ in \eqref{eq:input_FOH}, and then this computation is restarted at each discrete temporal node ${\bf x}_k$. 

The original optimal control problem in \ref{subsec:cont_time_ocp} is a free-final-time problem, i.e., the time of flight $t_{f}$ is a decision variable and not a fixed constant. To convert it into an equivalent fixed-final-time problem, we introduce a time scaling variable $\tau^k\in[0,1]$ that has a bijective mapping with the original time variable $t\in [t_k, t_{k+1}]$ as follows: 
\begin{align}
t = s_k \tau^k + t_k, \label{eq:time_scaling}
\end{align}
where $s_k$ is referred to as a time-interval dilation variable that is the length of the subinterval $s_k \coloneqq t_{k+1} - t_k$. With the conversion \eqref{eq:time_scaling}, the time-scaling variable $\tau^k$ represents the normalized time in $[0, 1]$ for each subinterval $[t_k,t_{k+1}]$. Henceforth, we use the simplified notation $\square(\tau^k)\coloneqq\square(s_k\tau^k+t_{k})$, where $\square$ is a placeholder for any time-varying variable, for brevity. Then, \eqref{eq:integration_equation} can be equivalently expressed with respect to $\tau^k\in[0,1]$ as
\begin{align}
{\bf x}(\tau^k) & ={\bf x}_{k}+\int_{0}^{\tau^k}s_k{\bf f}({\bf x}(s),{\bf u}(s))\text{d}s.\label{eq:integration_equation_tau}
\end{align}

\subsection{Constraints\label{subsec:Linearization}}

As discussed earlier, the proposed method linearizes all nonconvex constraints around the reference trajectory, which is either the solution to the subproblem at the previous iteration or the initial guess. The reference trajectory consists of the state profile $\bar{{\bf x}}_k$ and the control input profile $\bar{{\bf u}}_k$ for all $ k\in\mathcal{N}_{0}^{N}$, and the time-interval dilation variable $\bar{s}_k$ for all $ k\in\mathcal{N}_{0}^{N-1}$ that will be denoted together as $\bartraj$ for notational simplicity. Similarly, $\traj$ denotes the solution trajectory at the current iteration.

\subsubsection{Dynamics}
To ensure continuity of the state trajectory ${\bf x}(\tau^k)$ at each discrete-time node, we need to impose ${\bf x}_k = {\bf x}(t_k^-)$ where $t_k^-$ is the left limit of $t_k$ (See. Figure~\ref{fig:multiple_shooting}). This constraint can be written as:
\begin{align}
{\bf x}_{k+1} & ={\bf x}_{k}+\lim_{\tau\rightarrow 1^{-}}\int_{0}^{\tau}s_k{\bf f}({\bf x}(s),{\bf u}(s))\text{d}s, \label{eq:continuity_condition_tau}\end{align}
for all $k\in\mathcal{N}_{0}^{N}$ where ${\bf x}(\cdot)$ is the solution of \eqref{eq:integration_equation_tau}. With FOH interpolation \eqref{eq:input_FOH}, the right-hand side of \eqref{eq:continuity_condition_tau} is a nonconvex function of ${\bf x}_k, {\bf u}_k, {\bf u}_{k+1}, s_k$, and evaluating this function involves solving the ordinary differential equation (ODE). To impose the continuity constraint in our SCP framework, we need to linearize the constraint with respect to the reference trajectory $\bartraj$. One way to linearize \eqref{eq:continuity_condition_tau} is to use the variational method \cite{diehl2011numerical} that employs the sensitivity equation \cite{khalil2015nonlinear} given as the variational differential equation along the reference trajectory. More details on the derivation of the variational approach to inverse-free exact discretization can be found in \cite{kaya2003computational,loxton2008optimal,lin2014control,quirynen2015autogenerating,kamath2022real,kamath2023customized,mceowen2023high}.
Then, the linearization of the dynamics constraint 
 \eqref{eq:continuity_condition_tau} yields
\begin{align}
{\bf x}_{k+1} & =A_{k}{\bf x}_{k}+B_{k}^{m}{\bf u}_{k}+B_{k}^{p}{\bf u}_{k+1}+z_{k}s_{k}+w_{k},\label{eq:linearized_dynamic}
\end{align}
with $A_k \in \mathbb{R}^{n_x \times n_x}$, $B_k^m, B_k^p \in \mathbb{R}^{n_x \times n_u}$, $z_k \in \mathbb{R}^{n_x}$, and $w_k \in \mathbb{R}^{n_x}$.
\subsubsection{State and input constraints}
The state and input constraints given in \eqref{eq:const_height}-\eqref{eq:const_input_dot} are convex constraints, and as a result, it is not necessary to linearize them. These constraints can be imposed at each node point as follows:
\begin{subequations}
\label{eq:const_discrete_pointwise}
\begin{align}
    -(p_{D})_{k}&\geq0, \\
    {\bf v}_{\min}&\leq{\bf v}_{k}\leq{\bf v}_{\max},\quad \phi_{\min}\leq\phi_{k}\leq\phi_{\max}, \\
    \theta_{\min}&\leq\theta_{k}\leq\theta_{\max},\quad {\bf \Omega}_{\min}\leq{\bf \Omega}_{k}\leq{\bf \Omega}_{\max}, \\
    (u_{a})_{k}\tan\alpha_{\min}&\leq(w_{a})_{k}\leq(u_{a})_{k}\tan\alpha_{\max},\quad {\bf u}_{\min}\leq{\bf u}_{k}\leq{\bf u}_{\max}, \\
    \delta_{T,\min}&\leq(\delta_{T})_k\leq\delta_{T,\max},\quad\dot{\delta}_{\min}\leq\frac{1}{\tau_{T}}\left((\eta_{T})_{k}-(\delta_{T})_{k}\right)\leq\dot{\delta}_{\max}.
\end{align}
\end{subequations}

As mentioned earlier, the runway alignment constraint derived in \eqref{eq:runway_align_time} involves the existential quantifier for the switching time $t_s$. To incorporate this in the finite-dimensional optimization framework, we first specify a node point $k_s$ such that the time at the node is equal to the switching time, i.e., $t_s \defeq t_{k_s}$. Then, the runway alignment constraint can be represented as
\begin{align}
    -(p_{D})_k\leq h_{c} \text{ and }{\bf c}({\bf p}_k)\leq0, \quad \forall k\in \mathcal{N}_{k_s}^{N}.\label{eq:runway_align_node}
\end{align}
The alignment condition is imposed starting at the $k_s$-th node, and is imposed at all succeeding nodes. Note that although node point corresponding to $k_s$ is determined in advance, the switching time itself is not predetermined. Instead, the switching time is obtained by solving the optimization problem because the time-interval dilation variable $s_k$ is treated a decision variable. The switching time can be represented in terms of $t_{k_s}$ and $s_k$ as $t_{s} \defeq t_{k_s} = \sum_{k=0}^{k_s} s_k$. By imposing constraint \eqref{eq:runway_align_node}, we can guarantee that there exists a time after which the alignment condition is activated, thereby satisfying the original runway alignment constraint \eqref{eq:runway_align_time}. It is also worth noting that the constraint \eqref{eq:runway_align_node} is convex, and hence will not be linearized. 

Additionally, we utilize uniform meshes for each base leg phase and final approach phase. To be more specific, all nodes from the start to $k_s-1$ have identical time-interval dilation values $s_k$, that is,
\begin{align}
    s_i = s_j, \quad \forall i,j \in \mathcal{N}_0^{k_s-1}. \label{eq:constraint_on_s1}
\end{align}
This maintains a uniform mesh during the base leg phase. Similarly, all nodes from $k_s$ to $N-1$ also share the same $s_k$ values:
\begin{align}
    s_i = s_j, \quad \forall i,j \in \mathcal{N}_{k_s}^{N-1}. \label{eq:constraint_on_s2}
\end{align}
This ensures uniformity during the final approach phase. These constraints \eqref{eq:constraint_on_s1}-\eqref{eq:constraint_on_s2} may help to improve the convergence behavior by restricting the solution space.

Since the runway alignment constraint given by the STC in \eqref{eq:runway_align_STC} and the obstacle avoidance constraint in \eqref{eq:const_obstacle_avoidance} are nonconvex, they need to be linearized. The linear approximation of the runway alignment STC is given as follows:
\begin{align}
    {\bf g}({\bf p}) \approx
    \begin{cases}
        {\bf g}(\bar{\bf p}_k) + \frac{\partial {\bf g}}{\partial {\bf p}}\big|_{\bar{\bf p}_k}({\bf p}_k - \bar{\bf p}_k), & \text{if } -\bar{p}_D < \bar{h}_c,\\
        0, & \text{otherwise},
    \end{cases} \qquad \forall k \in \mathcal{N}_0^{N} \label{eq:STC_linear}
\end{align}
where $\bar{{\bf p}}_k \in \mathbb{R}^{3}$ is the position component of the reference state vector $\bar{{\bf x}}_k$. The obstacle avoidance constraint is linearized as follows:
\begin{align}
1-\norm{H_{i}(\bar{\bf p}_k-r_{i})}_{2}-\frac{H_{i}^{\top}H_{i}(\bar{\bf p}_k-r_{i})}{\norm{H_{i}(\bar{\bf p}_k-r_{i})}_{2}}({\bf p}_k-\bar{\bf p}_k)\leq0,\quad\forall i\in\mathcal{N}_{1}^{n_{\mathrm{obs}}}, \quad \forall k \in \mathcal{N}_0^{N}.\label{eq:obstacle_avoidance_linear}
\end{align}

\subsubsection{Continuous-time constraints} \label{subsub:CCS}
Enforcing constraints only at the node points does not guarantee their satisfaction in the subinterval between the node points, thereby leading to the inter-sample constraint violation \cite{dueri2017trajectory}. Here, we employ the method introduced in \cite{ctcs2024} and applied in \cite{nmpc2024,chari2024fast} to ensure continuous-time. Similar methods can also be found in \cite{teo1987simple,teo1989nonlinear,lin2014control}.

Let us first define a vector-valued function ${\bf h}$ representing constraints that need to be satisfied in continuous-time. Specifically, we focus on constraints related to roll and pitch angles. Our empirical observations indicate that violations of these constraints within a subinterval can lead to highly erratic and undesirable aircraft motion, unlike the other constraints in \eqref{eq:const_discrete_pointwise}. We define the function ${\bf h}:\mathbb{R}^{n_x}\rightarrow\mathbb{R}^{n_h}$ as
\begin{align}
    {\bf h}({\bf x}(t))&\coloneqq\left[\begin{array}{c}
    \phi - \phi_{\min}\\
    \phi_{\max} - \phi\\
    \theta - \theta_{\min}\\
    \theta_{\max} - \theta\\
    \end{array}\right],\label{eq:def_h}
\end{align}
with $n_h = 4$. Given the continuously differentiable function ${\bf h}$ and the continuous-time state solution ${\bf x}(\cdot)$, the following two expressions are equivalent \cite{ctcs2024}:
\begin{align}
    {\bf h}({\bf x}(t)) &\leq 0, \quad \forall t\in[t_k, t_{k+1}]
    \Leftrightarrow
    \int_{t_k}^{t_{k+1}} \max\{{\bf h}({\bf x}(t)),0\}^{2} \text{d}t = 0. \label{eq:const_on_h}
\end{align}
This implies that the pointwise inequality constraint in terms of time $t$ can be equivalently written as a single equality constraint with the integral representation. Using this equivalence, we impose the following constraint:
\begin{align}
   \int_{t_k}^{t_{k+1}} \max\{{\bf h}({\bf x}(t)),0\}^{2} \text{d}t = 0, \quad \forall k \in \mathcal{N}_0^{N-1}.\label{eq:CCT}
\end{align}
The nonconvex constraint \eqref{eq:CCT} can be linearized via the variational method, as applied to \eqref{eq:linearized_dynamic}. After linearization, we obtain
\begin{align}
\bar{A}_{k}{\bf x}_{k}+\bar{B}_{k}^{m}{\bf u}_{k}+\bar{B}_{k}^{p}{\bf u}_{k+1}+\bar{z}_{k}s_{k}+\bar{w}_{k}=0,  \quad \forall k\in\mathcal{N}_{0}^{N-1},\label{eq:CCT_linearized}
\end{align}
where $\bar{A}_k \in \mathbb{R}^{n_x \times n_x}$, $\bar{B}_k^m, \bar{B}_k^p \in \mathbb{R}^{n_x \times n_u}$, $\bar{z}_k \in \mathbb{R}^n_x$, and $\bar{w}_k \in \mathbb{R}^n_x$.

\subsection{Virtual control and trust region}
The linearization process described in \ref{subsec:Linearization} may cause: (1) artificial infeasibility, wherein the convex subproblem is rendered infeasible even if the original nonconvex problem is feasible, and (2) artificial unboundedness, wherein the subproblem cost function is unbounded from below, i.e., it can be minimized indefinitely. To address these issues, we introduce virtual control variables and impose a trust region penalty, respectively \cite{szmuk2020successive,reynolds2020dual,malyuta2022convex}. 

The virtual controls ${\bf \nu}^d_{k}\in\mathbb{R}^{n_{x}}$ and ${\bf \nu}^c_{k}\in\mathbb{R}^{n_{h}}$ are additional terms added to the linearized dynamics constraint \eqref{eq:linearized_dynamic} and the linearized continuous time constraints \eqref{eq:CCT_linearized}, respectively, as follows:
\begin{align}
{\bf x}_{k+1}&=A_{k}{\bf x}_{k}+B_{k}^{m}{\bf u}_{k}+B_{k}^{p}{\bf u}_{k+1}+z_{k}s+w_{k}+{\bf \nu}^d_{k},\quad \forall k\in\mathcal{N}_{0}^{N-1},\label{eq:linearized_dyn_vc} \\
0 &= \bar{A}_{k}{\bf x}_{k}+\bar{B}_{k}^{m}{\bf u}_{k}+\bar{B}_{k}^{p}{\bf u}_{k+1}+\bar{z}_{k}s_{k}+\bar{w}_{k} + {\bf \nu}^c_{k}\quad \forall k\in\mathcal{N}_{0}^{N-1}. \label{eq:linearized_cct_vc}
\end{align}
This modification with the virtual control terms ensures that the subproblem is always feasible. However, for the solution to be feasible with respect to the original nonlinear dynamics, the value of these terms should go to zero as the solution converges. To encourage this behavior, we impose the following penalty terms:
\begin{align}
J_{vc}=\sum_{k=0}^{N}\norm{{\bf \nu}^d_{k}}_{1} + \norm{{\bf \nu}^c_{k}}_{1}. \label{eq:cost_discrete_virtual}
\end{align}
The trust region penalty is set as follows:
\begin{align}
J_{tr}=\sum_{k=0}^{N}\norm{{\bf x}_{k}-\bar{{\bf x}}_{k}}_{2}^{2}+\norm{{\bf u}_{k}-\bar{{\bf u}}_{k}}_{2}^{2} + \abs{s_k - \bar{s}_k}^2. \label{eq:cost_discrete_trust}
\end{align}
The trust region penalty plays a role in keeping the solution in the vicinity of the reference trajectory around which the system is linearized, to ensure that the subproblem cost is not unbounded. Then, the conditions $J_{vc}\approx0$ and $J_{tr}\approx0$ are included in the stopping criteria, which will be discussed in the following subsection.

\subsection{Objectives in convex subproblem}

Here we illustrate how the cost functions defined in continuous-time domain, \eqref{eq:cost_all_cont}, are approximated to obtain a discretized convex subproblem. The cost functions for the final time and the thrust command in \eqref{eq:cost_time_thrust_cont} are defined as follows:
\begin{align}
J_{\text{time}}^{d} & =\sum_{k=0}^{N-1} s_k,\quad J_{\mathrm{thr}}^{d}=\sum_{k=0}^{N}(\eta_{T})^2_{k}, \label{eq:cost_discrete_final_input}
\end{align}
where we add the superscript ``$d$'' in order to distinguish from their continuous-time counterparts. Next, we consider the cost function for the rate of the control input commands in \eqref{eq:cost_time_control_rate_cont}. Since the control input command rates $\dot{\delta}_{A}$, $\dot{\delta}_{E}$, and $\dot{\delta}_{R}$ are not explicitly treated as control inputs, we penalize the change of the input commands, $\delta_{A}$, $\delta_{E}$, and $\delta_{R}$, as follows:
\begin{align}
J_{r}^{d}=\sum_{k=0}^{N-1}\left((\delta_{A})_{k+1}-(\delta_{A})_{k}\right)^{2}+\left((\delta_{E})_{k+1}-(\delta_{E})_{k}\right)^{2}+\left((\delta_{R})_{k+1}-(\delta_{R})_{k}\right)^{2}. \label{eq:cost_discrete_input_dot}
\end{align}
Finally, the cost function for the angular velocity is imposed as
\begin{align}
J_{\Omega}^{d}=\sum_{k=0}^{N}p_{k}^{2}+q_{k}^{2}+r_{k}^{2}. \label{eq:cost_discrete_pqr}
\end{align}
The overall objective function for the convex subproblem, based on \eqref{eq:cost_discrete_virtual}-\eqref{eq:cost_discrete_pqr}, is formulated as
\begin{align}
J^d=w_{t}J_{\text{time}}^{d}+w_{\mathrm{thr}}J_{\mathrm{thr}}^{d}+w_{r}J_{r}^{d}+w_{\Omega}J_{\Omega}^{d}+w_{vc}J_{vc}+w_{tr}J_{tr}, \label{eq:cost_discrete_total}
\end{align}
where $w_{vc}\in\mathbb{R}_{++}$ and $w_{tr}\in\mathbb{R}_{++}$ are weights for the virtual control and trust region penalties, respectively.

\begin{mybox}
    \vspace{-0.75em}
    \subsection{The discrete-time convex subproblem.} \label{subsec:summary_discrete_OCP}
    \vspace{-2.5em}
    \begin{align*}
    \underset{{\bf x}_{k},{\bf u}_{k},s_k,\nu^d_k,\nu^c_k}{\text{minimize}} \quad& \text{Eq \eqref{eq:cost_discrete_total}} \quad \text{subject to} \\
     \fbox{\text{Linearized dynamics}} \quad& \text{Eq \eqref{eq:linearized_dyn_vc}}, \\
     \fbox{\text{State and input constraints}}\quad& \text{Eqs \eqref{eq:const_discrete_pointwise}}, \\
     \fbox{\text{Runway alignment}} \quad& \text{Eq \eqref{eq:runway_align_node} or \eqref{eq:STC_linear}}, \\
     \fbox{\text{Constraints on mesh}} \quad& \text{Eqs \eqref{eq:constraint_on_s1}-\eqref{eq:constraint_on_s2}}, \\
     \fbox{\text{Obstacle avoidance}}\quad& \text{Eq \eqref{eq:obstacle_avoidance_linear}}, \\
     \fbox{\text{Continuous-time constraints}}\quad& \text{Eq \eqref{eq:linearized_cct_vc}}, \\
    \fbox{\text{Boundary conditions}} \quad& \text{Eqs \eqref{eq:boundary_initial}-\eqref{eq:boundary_final}}.
    \end{align*}
\end{mybox}

\subsection{Extrapolation update}
\label{subsec:extrapolation}
Let $\hattraj$ be the solution obtained by solving the convex subproblem in Section \ref{subsec:summary_discrete_OCP}. Then, the extrapolation step updates the solution as follows:
\begin{subequations}
\label{eq:extrapolation}
\begin{align}
    {\bf x}_k &= \bar{\bf x}_k + \gamma (\hat{\bf x}_k - \bar{\bf x}_k),\quad \forall k\in\mathcal{N}_{0}^{N}, \\
    {\bf u}_k &= \bar{\bf u}_k + \gamma (\hat{\bf u}_k - \bar{\bf u}_k),\quad \forall k\in\mathcal{N}_{0}^{N},  \\
    {s}_k &= \bar{s}_k + \gamma (\hat{s}_k - \bar{s}_k), \quad \forall k\in\mathcal{N}_{0}^{N-1}, 
\end{align}
\end{subequations}
where $\gamma \geq 1$ is an extrapolation parameter, and setting $\gamma = 1$ is equivalent to not performing extrapolation. The motivation behind this extrapolation is to exploit the direction provided by each subproblem by computing the difference between the subproblem solution and the reference solution, rather than simply accepting the subproblem solution itself. We aim to accelerate the convergence of solution by choosing $\gamma \geq 1$, motivated from Nesterov's acceleration technique \cite{nesterov1983method} and extrapolation-based first-order methods \cite{yang2017proximal,yu2022extrapolated}. This approach has been empirically validated through numerical simulations by comparing different values of $\gamma$.


\change{\subsection{Convergence analysis}}

\change{Without the extrapolation step ($\gamma = 1$), the solution method aligns with the framework of the prox-linear method \cite{drusvyatskiy2018error,ctcs2024}, enabling the exploitation of its convergence guarantees. The prox-linear method is designed to solve
\begin{align*}
    \Theta({\bf y}) = J({\bf y}) + H(G({\bf y})),
\end{align*}
where ${\bf y}$ represents the decision variable. Here, $H$ is a convex, $\alpha_h$-Lipschitz continuous function, and $G$ is potentially nonconvex continuously differentiable with a $\beta_g$-Lipschitz continuous gradient.  To optimize $\Theta({\bf y})$, the prox-linear method iteratively minimizes a convex approximation of $\Theta({\bf y})$, given by
\begin{align*}
    \Theta^\rho({\bf y};{\bf y}_i) = J({\bf y}) + H( G({\bf y}_i) + \nabla G({\bf y}_i) ({\bf y} - {\bf y}_i)) + \frac{1}{2\rho} \norm{{\bf y} - {\bf y}_i}_2^2,
\end{align*}
where the subscript $i$ denotes the current iteration index. In our specific setup, ${\bf y}$ is a concatenated vector of all decision variables, expressed as
\begin{align*}
    {\bf y} = ({\bf x}_1, {\bf u}_1, s_1, \ldots, {\bf x}_{N-1}, {\bf u}_{N-1}, s_{N-1}, {\bf x}_N, {\bf u}_N),
\end{align*}
and the constant $\rho$ is given by $1/(2w_{tr})$. The terms $J({\bf y})$, $H(G({\bf y}))$, and $G({\bf y})$ are 
\begin{align*}
    J({\bf y}) &= w_{t}J_{\text{time}}^{d}+w_{\mathrm{thr}}J_{\mathrm{thr}}^{d}+w_{r}J_{r}^{d}+w_{\Omega}J_{\Omega}^{d} + I_{\mathcal{Y}}({\bf y}),\\
    H(G({\bf y})) &= w_{vc}J_{vc},\\
    G({\bf y}) &= ({\bf x}_{1} - {\bf x}_{0} - \lim_{\tau\rightarrow 1}\int_{0}^{\tau}s_1{\bf f}({\bf x}(s),{\bf u}(s))\text{d}s, \cdots,
    {\bf x}_{N} - {\bf x}_{N-1} - \lim_{\tau\rightarrow 1}\int_{0}^{\tau}s_{N-1}{\bf f}({\bf x}(s),{\bf u}(s))\text{d}s, \\
     & \quad \ \ \int_{t_0}^{t_{1}} \max\{{\bf h}({\bf x}(t)),0\}^{2} \text{d}t, \cdots,
     \quad \ \ \int_{t_{N-1}}^{t_{N}} \max\{{\bf h}({\bf x}(t)),0\}^{2} \text{d}t),
\end{align*}
where $I_{\mathcal{Y}}$ is an indicator function, and $\mathcal{Y}$ represents a feasible set to convex constraints \eqref{eq:const_discrete_pointwise}-\eqref{eq:constraint_on_s2}. The prox-linear method has the following convergence result.
\begin{theorem}(Lemma 5.1 of \cite{drusvyatskiy2018error}) At iteration $i$ of the prox-linear method, the following holds
\begin{align*}
    \Theta ({\bf y}_i) \geq \Theta ({\bf y}_{i+1}) + \frac{\rho}{2}(2-\alpha_h\beta_g\rho) \norm{\mathcal{G}_\rho ({\bf y}_i)}^2,
\end{align*}
where $\mathcal{G}_\rho (\cdot)$ is the prox-gradient mapping defined as
\begin{align*}
    \mathcal{G}_\rho({\bf y}) = \frac{1}{\rho} \left(z - \operatorname*{argmin}_{{\bf y}'} \Theta^\rho ({\bf y}' ; {\bf y})\right).
\end{align*}
Hence, $\Theta^\rho({\bf y})$ is monotonically decreasing for $\rho \leq \frac{1}{\alpha _h\beta_g}$.
\end{theorem}
}

\change{With the extrapolation step ($\gamma > 1$), the convergence analysis of the method has not been investigated yet. However, we provide empirical validation of the benefits of using extrapolation through numerical simulations including Monte Carlo studies in Section \ref{sec:simulation}. In summary, the proposed method with the extrapolation step successfully reduces the number of iterations required for convergence by approximately 16.5$\%$. To guarantee convergence while preserving the practical benefits of extrapolation, we can switch the value of $\gamma$ to 1 after a few extrapolated iterations. Then, the final trajectory from the extrapolated algorithm could be viewed as a very good initial guess to the convergence-guaranteed prox-linear algorithm.
}


\subsection{Initial guess and algorithm summary}

\subsubsection{Initial trajectory guess}
\label{subsec:initialguess}
The initial trajectory guess is a user-specified trajectory—that need not be dynamically feasible—that is used as a reference trajectory at the first iteration. The typical approach is to use a straight-line guess obtained by linearly interpolating between the initial and final boundary points \cite{malyuta2022convex,szmuk2020successive}. Here, we employ the straight-line interpolation approach to construct the initial guess for the state and control 
input trajectory $\bar{\bf x}_{k}$, $\bar{\bf u}_{k}$ as follows:
\begin{subequations}
\label{eq:straightline}
\begin{align}
\bar{{\bf x}}_{k} & =\frac{N-k}{N}\bar{{\bf x}}_{i}+\frac{k}{N}\bar{{\bf x}}_{f},\\
\bar{{\bf u}}_{k} & =\bar{{\bf u}}_{i},\quad\forall k\in\mathcal{N}_{0}^{N},
\end{align}
\end{subequations}
where the vectors $\bar{\bf x}_i$ and $\bar{\bf x}_f$ are $({\bf p}_i,{\bf v}_i,{\bf \Phi}_i,{\bf \Omega}_i)$ and $(0,{\bf v}_f,{\bf \Phi}_f,0)$, respectively, and $\bar{{\bf u}}_{i}\in\mathbb{R}^{4}$ is the constant control input vector.
Additionally, it is also necessary to set the initial guess for the time-interval dilation variable $\bar{s}_k$. For instance, one can choose the following value:
\begin{align}
\bar{s}_k &= \frac{\bar{s}_{T}}{N-1}, \quad \bar{s}_T=\frac{\norm{{\bf p}_{i}}_{2}}{2\norm{{\bf v}_{\min}}_{2}}+\frac{\norm{{\bf p}_{i}}_{2}}{2\norm{{\bf v}_{\max}}_{2}}, \label{eq:initialguess_on_s}
\end{align}
where $\bar{s}_T$ represents the average of the flight-times considering minimum speed and maximum speed in a straight line between the initial point and the final point. 


\subsubsection{Algorithm summary}

The proposed SCP method for aircraft approach and landing is summarized in Algorithm \ref{alg:alg_summary}. The algorithm starts with the generation of an initial guess trajectory $\bartraj$, as discussed in Section \ref{subsec:initialguess}. Then linearization and discretization are performed, as illustrated in Section \ref{subsec:Linearization}, after which we obtain $\{A_{k}, B_{k}^{m}, B_{k}^{p}, z_{k}, w_k\}$ and $\{\bar{A}_{k}, \bar{B}_{k}^{m}, \bar{B}_{k}^{p}, \bar{z}_{k}, \bar{w}_k\}$. The next step is to solve the convex subproblem that is given in Section \ref{subsec:summary_discrete_OCP} to obtain $\hattraj$. The extrapolation step described in Section \ref{subsec:extrapolation} computes the updated solution $\traj$. The termination criteria to determine convergence of the solution is set as $\{(J_{vc}<\epsilon_{vc}) \text{ and } (J_{xtr}<\epsilon_{tr})\}$, where $\epsilon_{vc}, \epsilon_{tr}\in\mathbb{R}_{++}$ are user-specified tolerances for the virtual controls and the trust region, respectively. This procedure is repeated until the termination criterion is satisfied or the iteration count reaches the maximum defined value $I_{\max}\in\mathbb{R}_{++}$.

\begin{algorithm}
\begin{algorithmic}[1]
\State Initialize $\bartraj$ by \eqref{eq:straightline}-\eqref{eq:initialguess_on_s}
\For{$i=1,\ldots,I_{\max}$}
\State Compute $\{A_{k}, B_{k}^{m}, B_{k}^{p}, z_{k}, w_k\}$ and $\{\bar{A}_{k}, \bar{B}_{k}^{m}, \bar{B}_{k}^{p}, \bar{z}_{k}, \bar{w}_k\}$
\State Compute $\hattraj$ by solving the convex subproblem in Section \ref{subsec:summary_discrete_OCP} 
\State Compute $\traj$ via the extrapolation step in \eqref{eq:extrapolation}
\If{$J_{vc} < \epsilon_{vc}$ and $J_{xtr} < \epsilon_{tr}$}      
\State break
\EndIf
\State Update $\bartraj \leftarrow \traj$
\EndFor \\
\Return{$\traj$}
\end{algorithmic}
\caption{Extrapolated PTR (xPTR)}
\label{alg:alg_summary}
\end{algorithm}

\section{Numerical Simulations} \label{sec:simulation}
In this section, we validate the proposed method through extensive numerical simulations and highlight the main contributions of the work. In Section \ref{subsec:vs_align}, we present trajectory results starting from 3 different initial conditions, and test the effectiveness of the runway alignment constraint by comparing the results simulated with and without the constraint. Also, we compare the performance of the algorithm and the resulting trajectories based on different implementations of the runway alignment constraint, one given by the proposed method \eqref{eq:runway_align_time} and the other, by the standard STC method \eqref{eq:runway_align_STC}.
The necessity of continuous-time constraint satisfaction is investigated in Section \ref{subsec:vs_CTC} by showing the case where enforcing the constraint only at node points leads to the large inter-sample constraint violation, whereas the proposed method leads to no such violation. In Section \ref{subsec:vs_CTC}, we also test an obstacle avoidance scenario, where there are two cylindrical obstacles between the initial and the final points. Next, we generate trajectories for the case in which the aircraft maneuvers under constant wind fields, and discuss how the winds affect the resulting trajectories in Section \ref{subsec:vs_wind}. Finally, a Monte Carlo study according to different initial conditions and wind fields are conducted in Section \ref{subsec:monte}, where we measure success rate, total number of iterations, and computation time.

The parameters for the aircraft dynamics and the SCP algorithm are given in Tables \ref{tab:dynamics} and \ref{tab:others}, respectively. 
All numerical simulations are run on an Apple MacBook Pro with an M1 Pro and an 8-core CPU. To solve the convex subproblem described in Section \ref{subsec:summary_discrete_OCP}, we use the Gurobi solver \cite{gurobi2018gurobi} using the CVXPY \cite{diamond2016cvxpy} parser interface in Python.

\begin{table}[H]
\caption{Simulation parameters of aircraft dynamics}
\label{tab:dynamics}
\centering{}%
\begin{tabular}{cccc}
\toprule 
Parameters & Value (Unit) & Parameters & Value (Unit)\tabularnewline
\midrule 
$m$ & $120000$ (kg) & g & $9.81$(m/s$^{2}$)\tabularnewline
\midrule 
$(J_{x},J_{y}J_{z},J_{xz})/m$ & $(40.07,64,99.92,2.0923)$ (m$^{2}$) & $\rho$ & $1.225$ (kg/m$^{3}$)\tabularnewline
\midrule 
$C_{L_{\alpha}},C_{L_{\alpha_{t}}},C_{L_{qV}}$ & $5.5,3.1,1.3$ ($\cdot$) & $\alpha_{0},\epsilon_{\alpha}$ & $-11.5,5.5$ (rad)\tabularnewline
\midrule 
$S_{t},S$ & $64,260$ (m$^{2}$) & $C_{D0},C_{D1},C_{D2},C_{D_{\alpha}}$ & 0.13, 0.654, 0.07, 5.5 ($\cdot$)\tabularnewline
\midrule 
$C_{Y\beta},C_{Y\delta_{R}}$ & $-1.6,0.24$ ($\cdot$) & $\bar{c}$ & $6.6$ (m)\tabularnewline
\midrule 
$C_{l\beta},C_{lp},C_{lr}$ & $-1.4,-11,5$ ($\cdot$) & $C_{l\delta_{A}},C_{l\delta_{R}}$ & $-0.6,0.22$ ($\cdot$)\tabularnewline
\midrule 
$C_{m0},C_{m\alpha},C_{mq},C_{m\delta_{E}}$ & $-0.59,-3.1,-4.03,-3.1$ ($\cdot$) & $C_{n\beta},C_{n\alpha\beta},C_{np},C_{nr}$ & $1,-3.82,1.7,-11.5$ ($\cdot$)\tabularnewline
\midrule 
$C_{n\delta_{A}},C_{n\delta_{R}}$ & $0,-0.63$ ($\cdot$) & ${\bf p}_{cg}$ & $(0.23\bar{c},0,0.1\bar{c})$ (m)\tabularnewline
\midrule 
${\bf p}_{ac}$ & $(0.12\bar{c},0,0)$ (m) & $\tau_{t}$ & $1.5$ (s)\tabularnewline
\bottomrule
\end{tabular}
\end{table}

\begin{table}[H]
\caption{Simulation parameters of constraints, objectives, and tolerances}
\label{tab:others}
\centering{}%
\begin{tabular}{cccc}
\toprule 
Parameters & Value (Unit) & Parameters & Value (Unit)\tabularnewline
\midrule 
${\bf v}_{\min}$ & $(80,-4,-10)$ (m/s) & ${\bf v}_{\max}$ & $(120,4,10)$ (m/s)\tabularnewline
\midrule 
$\phi_{\min},\phi_{\max}$ & $-15,15$ (deg) & $\theta_{\min},\theta_{\max}$ & $-15,15$ (deg)\tabularnewline
\midrule 
${\bf \Omega}_{\min}$ & $(-10,-10,-10)$ (deg/s) & ${\bf \Omega}_{\max}$ & $(10,10,10)$ (deg/s)\tabularnewline
\midrule 
$\alpha_{\min},\alpha_{\max}$ & $-11.5,14.5$ (deg) & $\delta_{T,\min},\delta_{T,\max}$ & $0,10$ (deg)\tabularnewline
\midrule 
${\bf u}_{\min}$ & $(-25,-25,-30,0)$ (deg) & ${\bf u}_{\max}$ & $(25,10,30,10)$ (deg)\tabularnewline
\midrule 
$\dot{\delta}_{\min},\dot{\delta}_{\max}$ & $-0.53,0.53$ (deg/s) & ${\bf v}_{f}$ & $(85,3,3)$ (m/s)\tabularnewline
\midrule 
$\theta_{\mathrm{lat}}^{\min},\theta_{\mathrm{lat}}^{\max}$ & $-2,2$ (deg) & $\theta_{\mathrm{ver}}^{\min},\theta_{\mathrm{ver}}^{\max}$ & $(3,5)$ (deg)\tabularnewline
\midrule 
$h_{c}$ & $500$ (m) & $w_{t},w_{\mathrm{thr}},w_{r},w_{\Omega}$ & $0.02,10,1,10$ ($\cdot$)\tabularnewline
\midrule 
$\epsilon_{vc}$ & $10^{-6}$ ($\cdot$) & $\epsilon_{tr}$ & $10^{-3}$ ($\cdot$)\tabularnewline
\bottomrule
\end{tabular}
\end{table}

\subsection{Runway alignment constraint} \label{subsec:vs_align}

\begin{figure}
    \centering
        \centering
         \begin{subfigure}[b]{0.4\textwidth}
         \centering
         \includegraphics[width=\textwidth,clip]{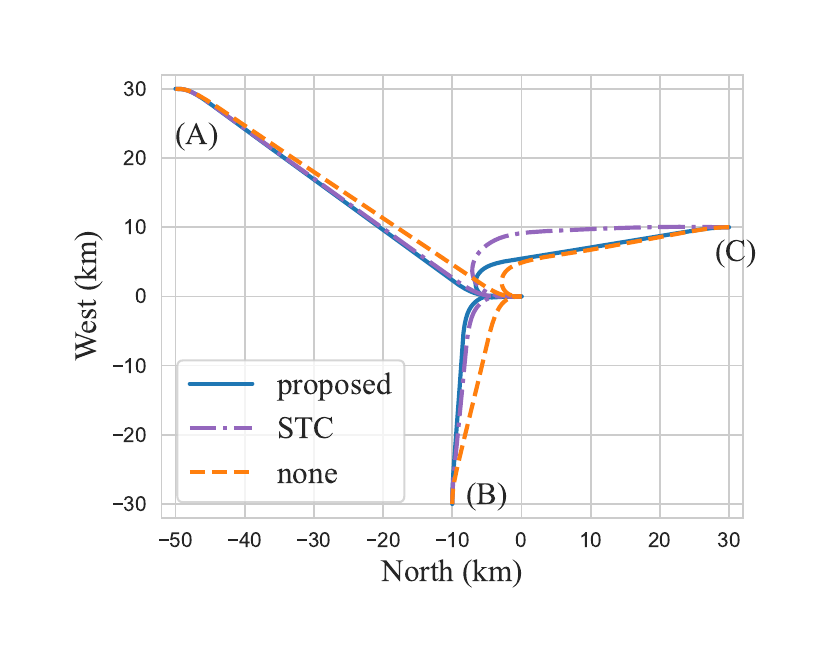}
         \end{subfigure}
         \begin{subfigure}[b]{0.5\textwidth}
         \centering
         \includegraphics[width=0.8\textwidth,trim={2cm 3cm 0.2cm 3.5cm},clip]{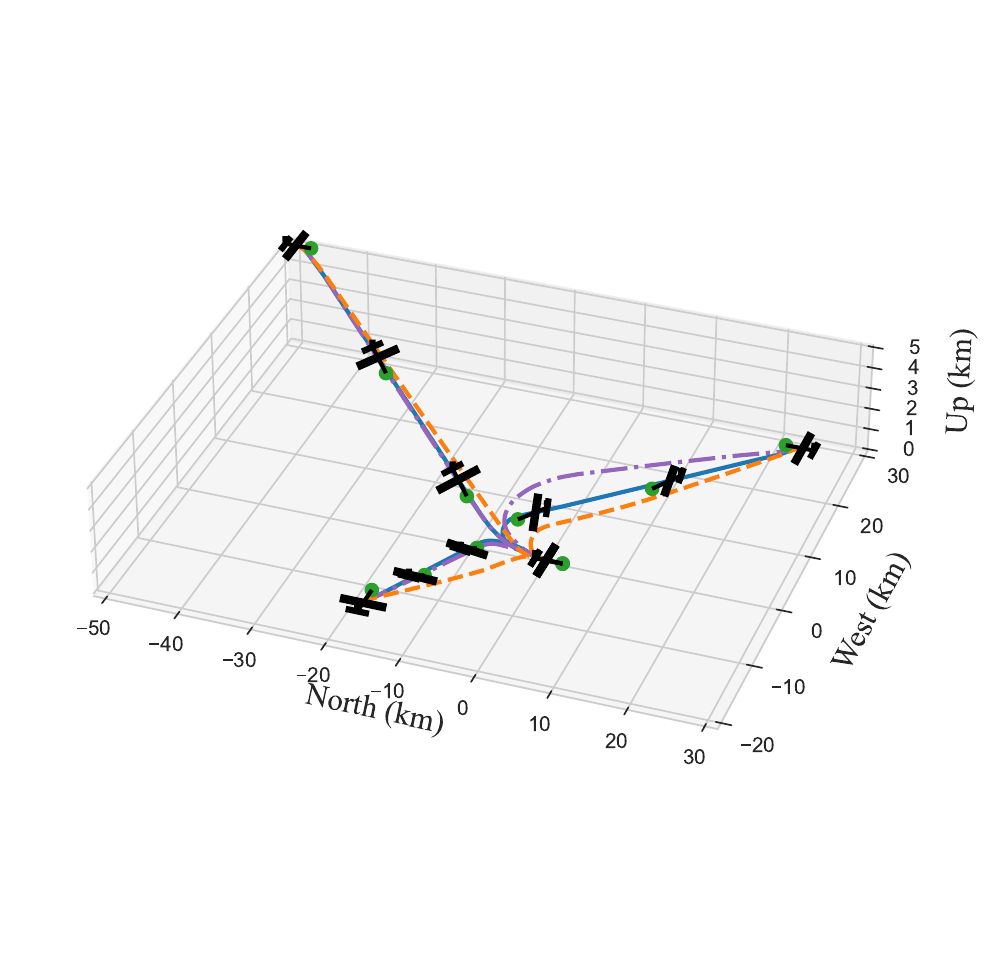}
         \end{subfigure}
\caption{Figure of trajectory results. Left: position results projected onto the North-West plane. Right: position results with respect to the North-West-Up frame.}
\label{fig:01_traj_result}
\end{figure}

\begin{figure}
    \centering
         \begin{subfigure}{1.0\textwidth}
         \centering
         \includegraphics[width=\textwidth,clip]{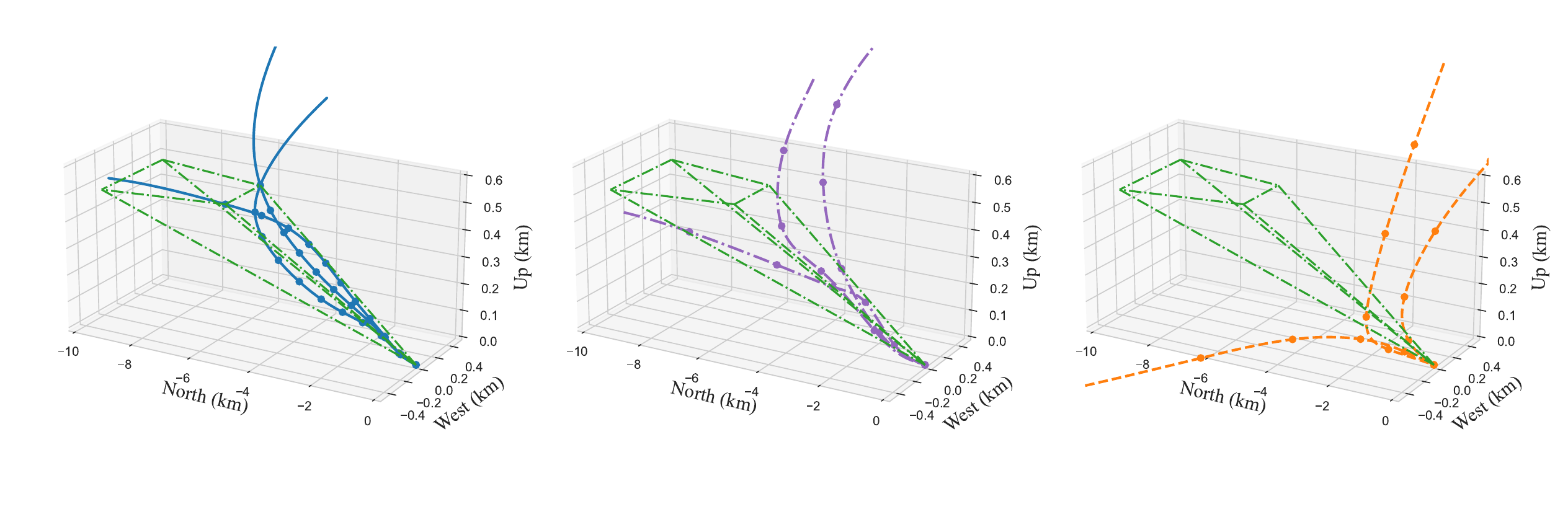}
         \end{subfigure}
         \begin{subfigure}{1.0\textwidth}
         \centering
         \includegraphics[width=\textwidth,clip]{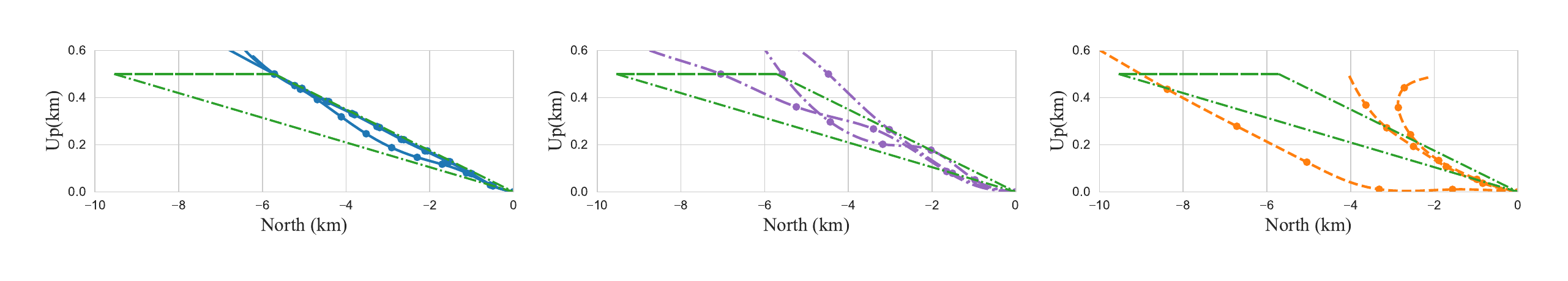}
         \end{subfigure}
         \begin{subfigure}{1.0\textwidth}
         \centering
         \includegraphics[width=\textwidth,clip]{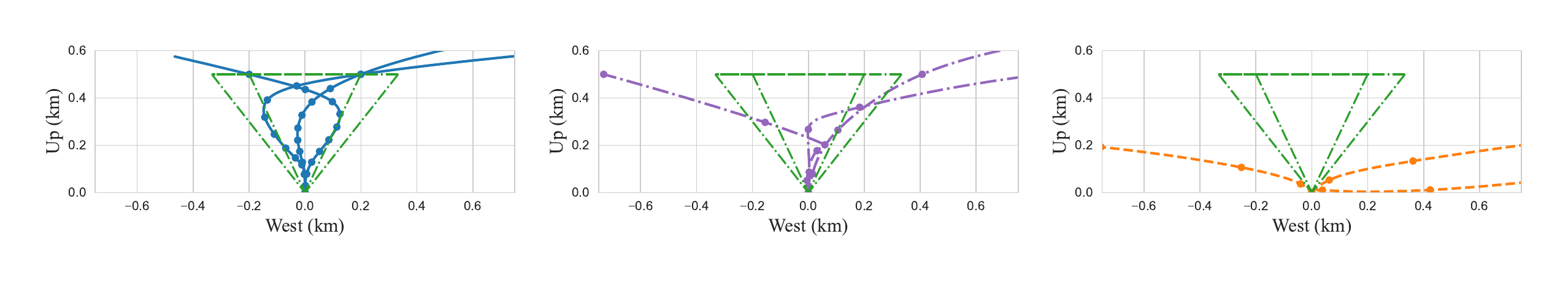}
         \end{subfigure}
\caption{Close-up views of the runway from Figure \ref{fig:01_traj_result}. Further details are provided in the main text.}
\label{fig:01_zoomin}
\end{figure}

\begin{figure}
    \centering
        \centering
        \includegraphics[width=0.9\textwidth,clip]{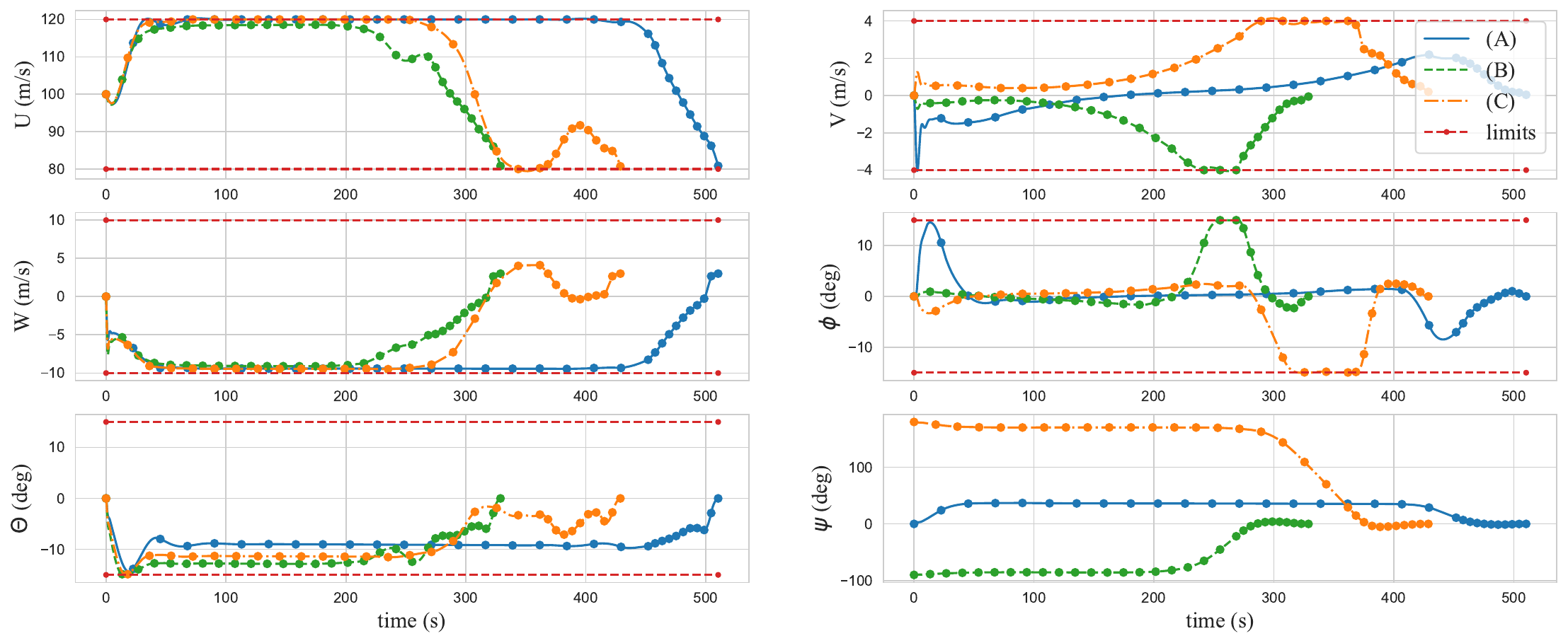}
\caption{Velocity and Euler angle profiles generated using the proposed method.}
\label{fig:01_vel_att}
\end{figure}

\begin{figure}
    \centering
        \centering
        \includegraphics[width=0.9\textwidth,clip]{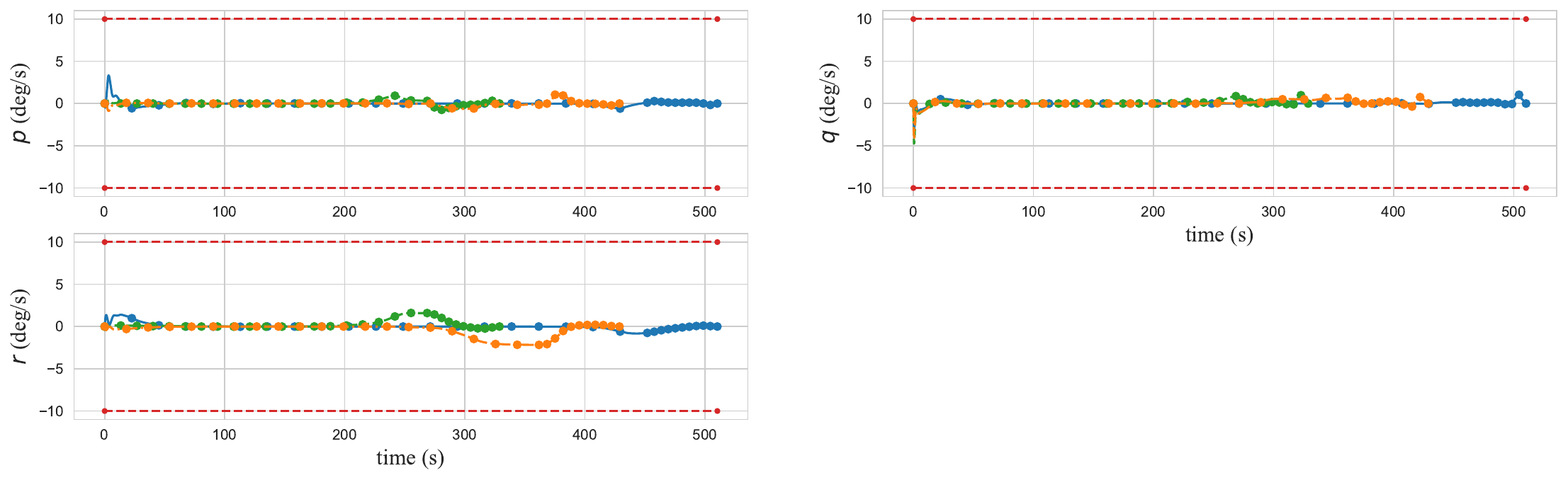}
\caption{Angular velocity profiles generated using the proposed method.}
\label{fig:01_ang_vel}
\end{figure}

\begin{figure}
    \centering
        \centering
        \includegraphics[width=0.9\textwidth,clip]{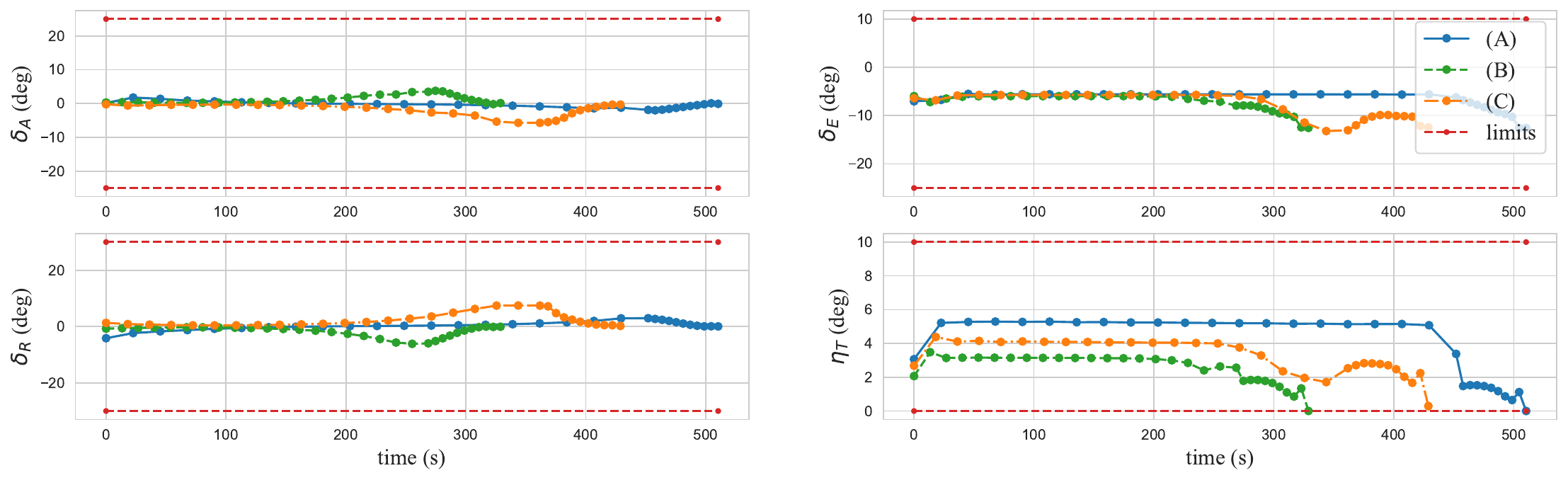}
\caption{Control input profiles generated using the proposed method.}
\label{fig:01_input}
\end{figure}

Figures \ref{fig:01_traj_result}-\ref{fig:01_input} illustrate the trajectory results computed from the following three different initial conditions:
\begin{subequations}
\label{eq:simulation_initial_conditions}
\begin{align}
    \text{(A)}: {\bf p}_i &= (-50,-30,-5)\ (\text{km}), {\bf v}_i=(100,0,0)\ (\text{m/s}), {\bf \Phi}_i = (0,0,0)\ (\text{deg}),{\bf \Omega}_i = 0_{3\times1}\ (\text{deg/s}),\\ 
    \text{(B)}: {\bf p}_i &= (-10,30,-5)\ (\text{km}),{\bf v}_i=(100,0,0)\ (\text{m/s}), {\bf \Phi}_i = (0,0,-90)\ (\text{deg}),{\bf \Omega}_i = 0_{3\times1}\ (\text{deg/s}),\\  
    \text{(C)}: {\bf p}_i &= (30,-10,-5)\ (\text{km}),{\bf v}_i=(100,0,0)\ (\text{m/s}), {\bf \Phi}_i = (0,0,90)\ (\text{deg}),{\bf \Omega}_i = 0_{3\times1}\ (\text{deg/s}).
\end{align}
\end{subequations}
The position results are given in Figure \ref{fig:01_traj_result}, with a detailed zoom-in of the runway provided in Figure \ref{fig:01_zoomin}. In Figure \ref{fig:01_zoomin}, the pentahedron represented by dash-dotted green lines is the feasible region with respect to the runway alignment constraint. The first column from the left shows the result of the proposed method, the second column shows that of the standard STC method, and the third column shows the case without the runway alignment constraint. For each case, the second and the third plots from the top show the position results projected onto North-Up and West-Up planes, respectively. Circular marks on trajectories illustrate node points.

The results of the total cost \eqref{eq:cost_discrete_total} and the iteration count are summarized in Table~\ref{tab:iteration_cost}. All trajectories are generated with the identical parameters and tolerances as detailed in Table~\ref{tab:others} with the weights $w_{vc}=10^2$ and $w_{tr}=1$ and the extrapolation parameter $\gamma = 1$, with the exception of the trajectory computed using the STC method starting point (B). 
The velocity and Euler angle profiles are shown in Figure \ref{fig:01_vel_att}, and the angular velocity and control input profiles are illustrated in Figure \ref{fig:01_ang_vel} and Figure \ref{fig:01_input}, respectively.

It is clear from Figure \ref{fig:01_traj_result} and Figure \ref{fig:01_zoomin} that the trajectories computed using the proposed method are aligned with the runway around the touch down point, as intended, and hence they are more operational. Meanwhile the trajectories computed without the runway alignment constraint tend to either be too steep or too shallow during the final landing phase, neither of which are reasonable in practice. The cost function values for the cases without the alignment constraint are smaller on average than those for the cases with the alignment constraint. This is expected since a trajectory with fewer constraints can be expected to have a lower cost. The trajectories computed using the standard STC method are aligned with the runway as well, but they tend to have larger iteration counts. \change{As pointed out previously, this is due to having the STC constraint implemented as a nonconvex constraint \eqref{eq:runway_align_STC}, thereby making the feasible set more challenging to optimize over. Furthermore, we can see in Figure \ref{fig:01_zoomin} that the trajectories computed using the standard STC method exhibit large inter-sample constraint violations with respect to the alignment constraint. One remedy for this is to combine the method for the continuous-time constraint satisfaction given in Section \ref{subsec:extrapolation} with the STC constraint, but this combination will further exacerbate the nonconvexity of the problem.}

\begin{table}[H]
\caption{Iteration count and cost results.}
\label{tab:iteration_cost}
\begin{centering}
\begin{tabular}{ccccc}
\toprule
\multicolumn{2}{c}{} & Proposed & STC & None\tabularnewline
\hline
\multirow{2}{*}{(A)} & Iteration & 26 & 28  & 23 \tabularnewline
\cline{2-5} \cline{3-5} \cline{4-5} \cline{5-5}
 & Cost & 11.9257  & 12.3158  & 12.1952 \tabularnewline
\hline
\multirow{2}{*}{(B)} & Iteration & 29 & 69  & 31 \tabularnewline
\cline{2-5} \cline{3-5} \cline{4-5} \cline{5-5}
 & Cost & 7.2956  & 7.2119  & 6.5909 \tabularnewline
\hline
\multirow{2}{*}{(C)} & Iteration & 26 & 67  & 28 \tabularnewline
\cline{2-5} \cline{3-5} \cline{4-5} \cline{5-5}
 & Cost & 9.7609  & 10.9931  & 8.1699 \tabularnewline
\hline
\end{tabular}
\par\end{centering}
\end{table}

\subsection{Continuous-time constraint satisfaction and obstacle avoidance} \label{subsec:vs_CTC}

\begin{figure}
    \centering
        \centering
         \begin{subfigure}[b]{0.4\textwidth}
         \centering
         \includegraphics[width=\textwidth,clip]{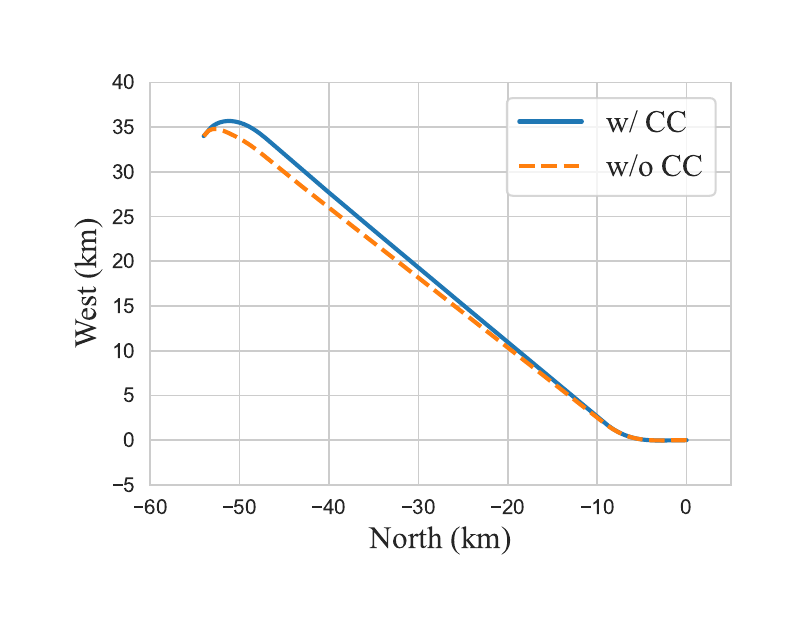}
         \end{subfigure}
         \begin{subfigure}[b]{0.5\textwidth}
         \centering
         \includegraphics[width=0.8\textwidth,trim={2cm 3cm 0.2cm 3.5cm},clip]{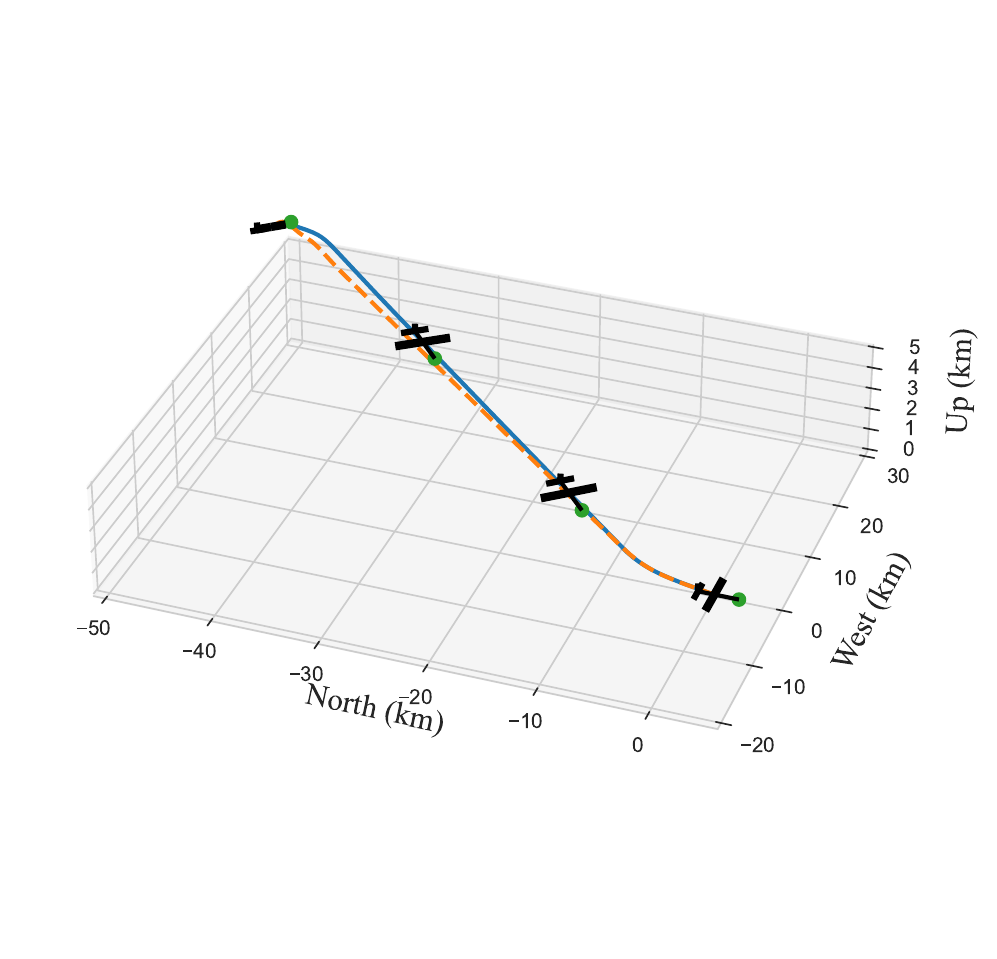}
         \end{subfigure}
\caption{Trajectory results illustrating the comparison between the cases with and without the continuous-time constraint (CC).}
\label{fig:02_traj_result}
\end{figure}

\begin{figure}
    \centering
        \centering
        \includegraphics[width=0.9\textwidth,clip]{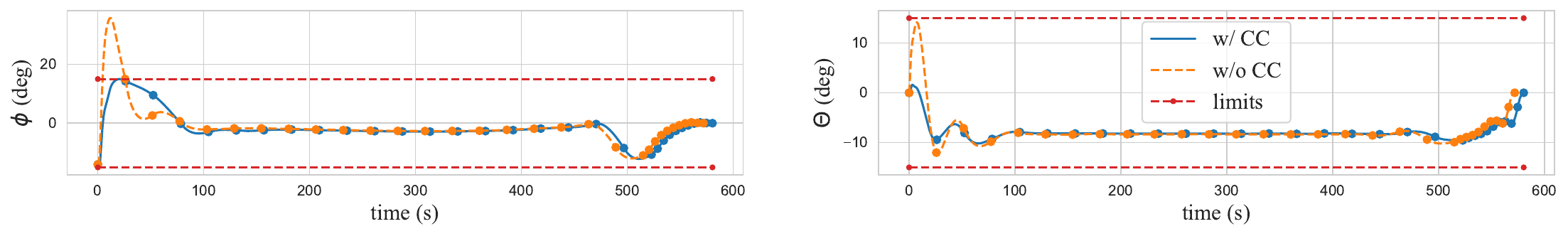}
\caption{The roll and pitch angle profiles with and without the continuous-time constraint (CC).}
\label{fig:02_att}
\end{figure}

We validate the proposed method in terms of continuous-time constraint satisfaction (see Section \ref{subsub:CCS}). To this end, we generate trajectories with and without the continuous-time constraints \eqref{eq:linearized_cct_vc}, starting from the following initial condition:
\begin{align}
    {\bf p}_i &= (-54,-34,-4.7)\ (\text{km}), {\bf v}_i=(100,0,0)\ (\text{m/s}), {\bf \Phi}_i = (-14,0,-50)\ (\text{deg}),{\bf \Omega} = 0_{3\times1}\ (\text{deg/s}).
\end{align}
The position results are shown in Figure \ref{fig:02_traj_result} and the roll and the pitch angle results are shown in Figure \ref{fig:02_att}. The trajectories are computed with the weights $w_{vc}=10^2$ and $w_{tr}=1$ and the extrapolation parameter $\gamma = 1$.

It is clear that the trajectory without the continuous-time constraint violates the roll limit severely during the first subinterval although the constraints are satisfied at each node point. The roll angle becomes even larger than 30 (deg), which is over twice the imposed limit. On the other hand, the result with the proposed method satisfies the roll and pitch angle limits during the entire flight. This shows the necessity of continuous-time constraint satisfaction and the effectiveness of the proposed method.

\begin{figure}
    \centering
        \centering
         \begin{subfigure}[b]{0.4\textwidth}
         \centering
         \includegraphics[width=\textwidth,clip]{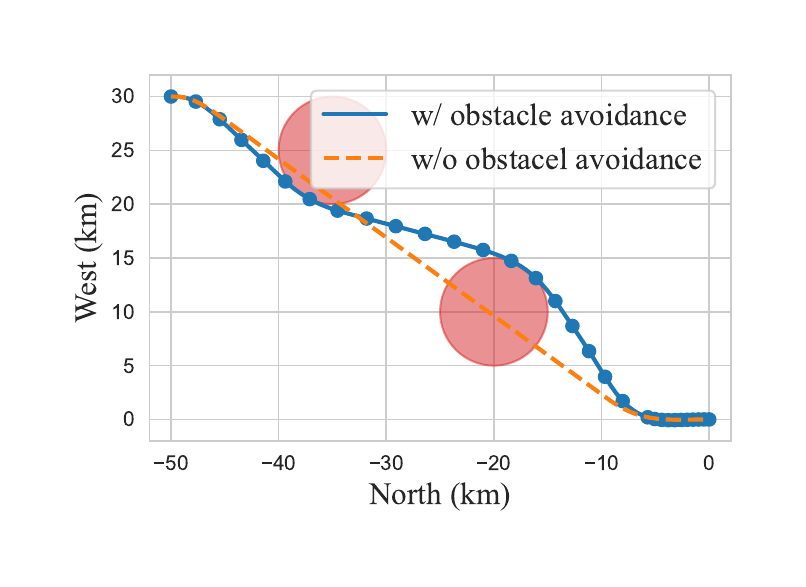}
         \end{subfigure}
         \begin{subfigure}[b]{0.5\textwidth}
         \centering
         \includegraphics[width=0.9\textwidth,trim={2cm 3cm 0.2cm 3.5cm},clip]{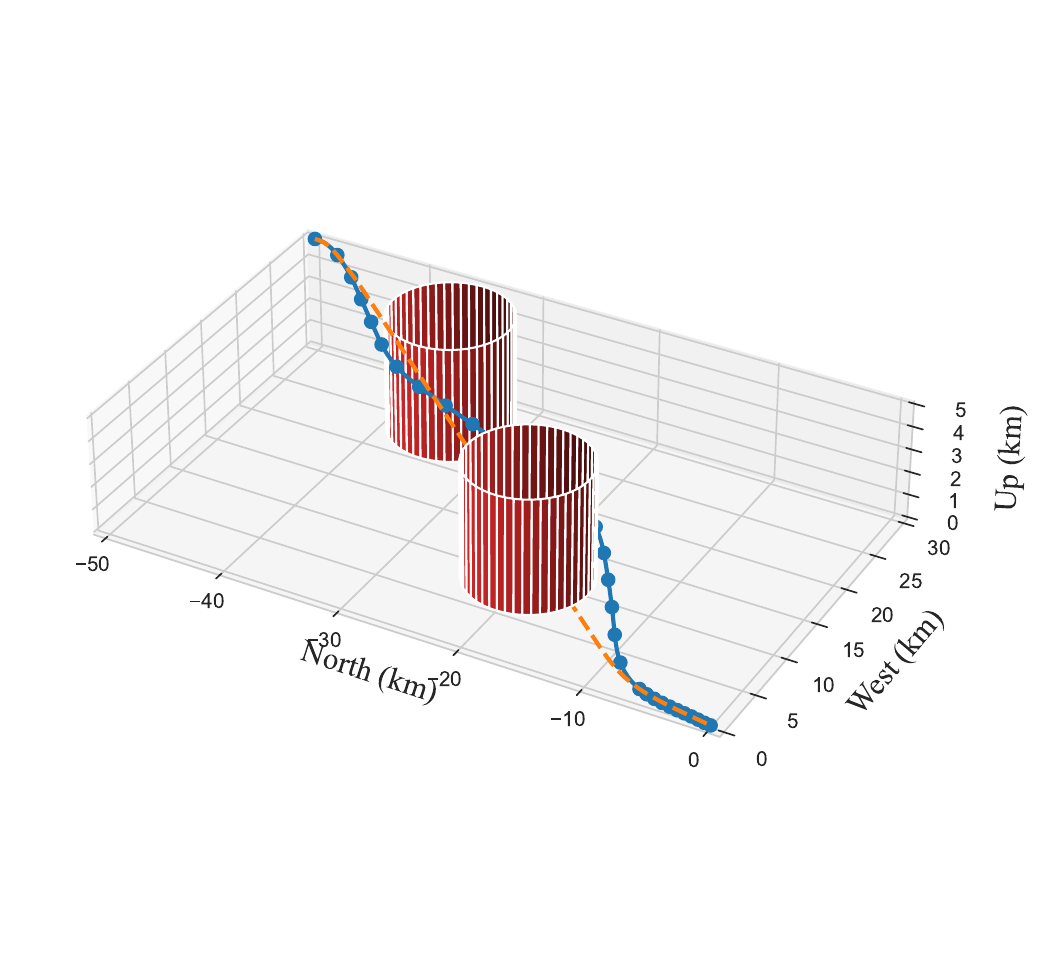}
         \end{subfigure}
\caption{Figure of trajectory results with obstacle avoidance.}
\label{fig:03_traj_result}
\end{figure}

Figure \ref{fig:03_traj_result} shows the trajectory result starting at the initial condition (A) given in \eqref{eq:simulation_initial_conditions} with the obstacle avoidance constraint \eqref{eq:const_obstacle_avoidance}. In the environment, there exist two cylindrical obstacles between the initial and the touch down points. The resulting trajectory successfully avoids these obstacles.

\subsection{Trajectory generation under constant wind fields}  \label{subsec:vs_wind}

\begin{figure}
    \centering
        \centering
         \begin{subfigure}[b]{0.4\textwidth}
         \centering
         \includegraphics[width=\textwidth,clip]{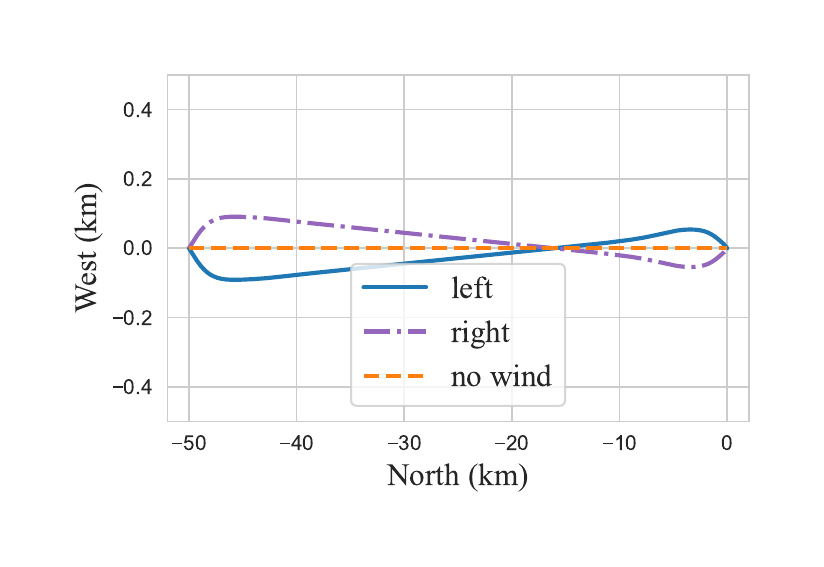}
         \end{subfigure}
         \begin{subfigure}[b]{0.5\textwidth}
         \centering
         \includegraphics[width=0.8\textwidth,trim={2cm 3cm 0.2cm 3.5cm},clip]{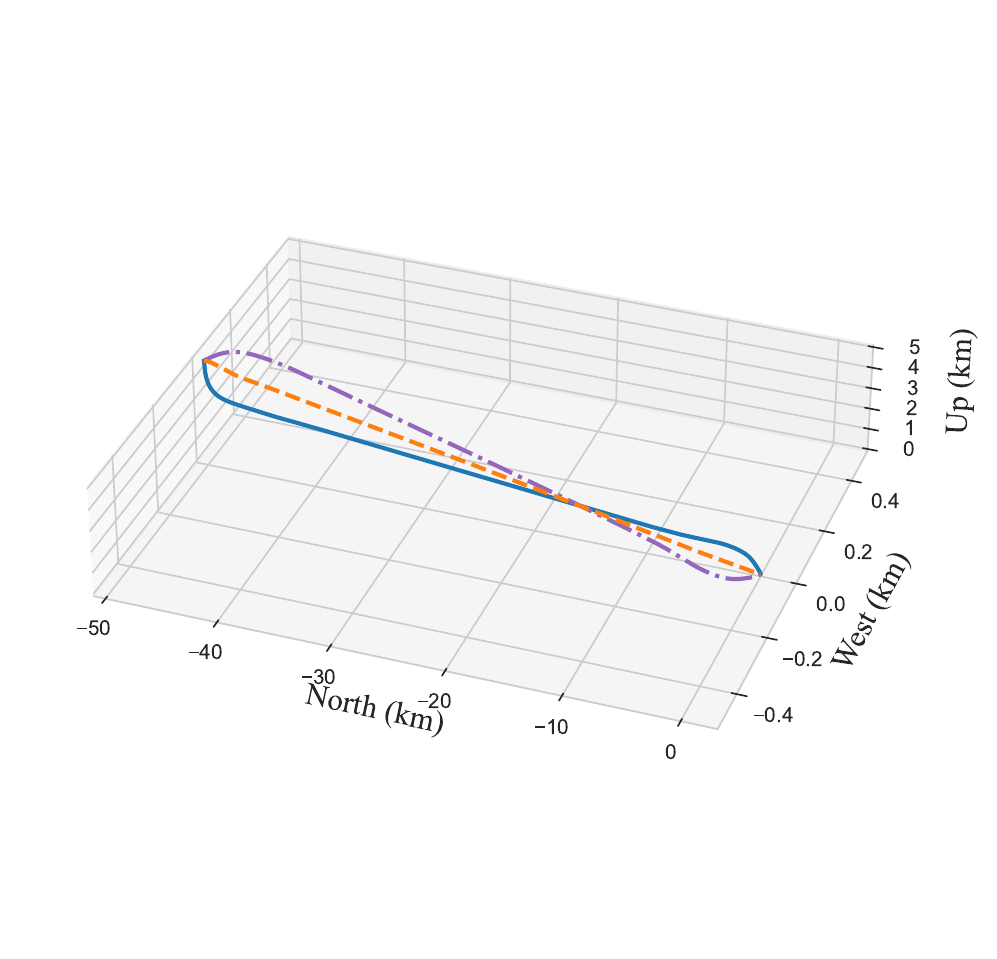}
         \end{subfigure}
\caption{Trajectory results with crosswind.}
\label{fig:04_traj_result}
\end{figure}

\begin{figure}
    \centering
        \centering
        \includegraphics[width=0.9\textwidth,clip]{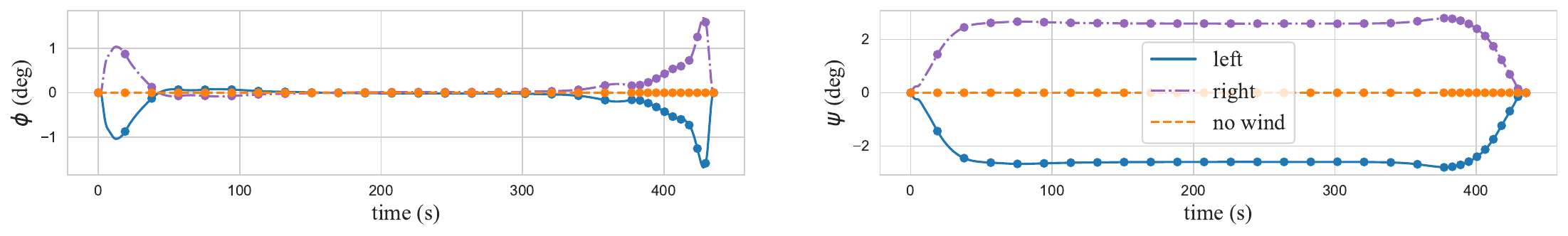}
\caption{The roll and yaw angle profiles with crosswind.}
\label{fig:04_att}
\end{figure}

In this case, we generate trajectories under crosswind conditions where the wind velocity vector is set to ${\bf w}=(0,5,0)$ (m/s) and ${\bf w} = (0,-5,0)$ (m/s) to represent the left (west-to-east) and right (east-to-west) crosswinds, respectively. All trajectories have the following initial condition:
\begin{align}
    {\bf p}_i &= (-50,0,-5)\ (\text{km}), {\bf v}_i=(100,0,0)\ (\text{m/s}), {\bf \Phi}_i = 0_{3\times1}\ (\text{deg}),{\bf \Omega}_i = 0_{3\times1}\ (\text{deg/s}).
\end{align}
The  position profile and roll and yaw angle profiles are illustrated in Figure \ref{fig:04_traj_result} and Figure \ref{fig:04_att}, respectively. These trajectories are computed with the weights $w_{vc}=10^2$ and $w_{tr}=1$ and the extrapolation parameter $\gamma = 1$.

To accomplish approach and landing in the presence of crosswinds, the resulting trajectory employs both the crab method and the wing-low method \cite{federal2011airplane}. For example, under the left crosswind condition, the aircraft initially shifts to the east-side slightly due to the left wind, and then maintains a negative yaw angle that makes a heading toward the wind so that the ground velocity ${\bf v}_g$ is aligned with and forwarded to the runway. This discrepancy between the aircraft's longitudinal direction and the direction of its ground velocity gives the appearance that the aircraft is performing a crabbing maneuver. Furthermore, we can see that once the aircraft gets closer to the runway during the final approach, the aircraft lowers the upwind wing (left wing), resulting in a negative roll angle. This wing-low method helps keep the aircraft's heading aligned with the runway under the crosswind condition, so that the aircraft can satisfy the zero yaw angle condition at touch down.

\subsection{Monte Carlo study} \label{subsec:monte}

In this subsection, we assess the robustness of the proposed algorithm by means of a Monte Carlo study. We consider a wide range of initial states and wind conditions, and also compare the performance of the algorithm based on different extrapolation parameters, denoted as $\gamma$.

To choose the initial states for the study, we perturb the position ${\bf p}$, the $x$-axis velocity $u$, and the Euler angles ${\bf \Phi}$. We generate 100 samples around each initial point (A), (B), and (C) specified in \eqref{eq:simulation_initial_conditions}, yielding a total of 300 initial states for our analysis. The position perturbation is allowed to take values in the range $\{(-5,-5,-0.5)\ \text{(km)}, (5,5,0.5)\ \text{(km)}\}$, the $x$-axis velocity perturbation is allowed to take values in the range $\{-10\ \text{(m/s)}, 10\ \text{(m/s)}\}$, and the attitude perturbation is allowed to take values in the range  $\{(-15,-15,-90)\ \text{(deg)}, (15,15,90)\ \text{(deg)}\}$. Separately, the wind vector ${\bf w}$ is randomly selected such that $\norm{{\bf w}}_{2} \leq 5\ \text{(m/s)}$.

\begin{algorithm}
\caption{Outcome Categorization Algorithm}
\begin{algorithmic}[1]
\State \textbf{Input:} initial condition (${\bf p}_{i}, {\bf v}_{i}, {\bf \Phi}_{i}, {\bf \Omega}_i$)
\State \textbf{Output:} Outcome label
\State \textbf{Initialize:} $w_{vc}=10^2$, $w_{tr}=1$, $I_{\max}=100$
\For{each initial condition (${\bf p}_{i}, {\bf v}_{i}, {\bf \Phi}_{i}, {\bf \Omega}_i$)}
    \For{attempt $a = 1$ to $2$}
        \State Run algorithm to check convergence for attempt $a$
        \If{attempt $a$ converges}
            \State Label as ``Success''
            \State \textbf{break}
        \ElsIf{attempt $a$ reaches $I_{\max}$}
            \If{$J_{vc} \geq \epsilon_{vc}$}
                \State Double $w_{vc}$
            \EndIf
            \If{$J_{tr} \geq \epsilon_{tr}$}
                \State Double $w_{tr}$
            \EndIf
            \If{$J_{vc} \geq \epsilon_{vc}$ \textbf{or} $J_{tr} \geq \epsilon_{tr}$ \textbf{and} $a == 2$}
                \State Label as ``Fail by max iteration''
            \EndIf
        \Else 
            \If{$a == 1$} 
                \State Double $w_{tr}$
            \Else
                \State Label as ``Fail by divergence''
            \EndIf
        \EndIf
    \EndFor
\EndFor
\end{algorithmic}
\label{alg:outcome_label}
\end{algorithm}

\begin{figure}
    \centering
        \centering
         \begin{subfigure}[b]{0.4\textwidth}
         \centering
         \includegraphics[width=\textwidth,clip]{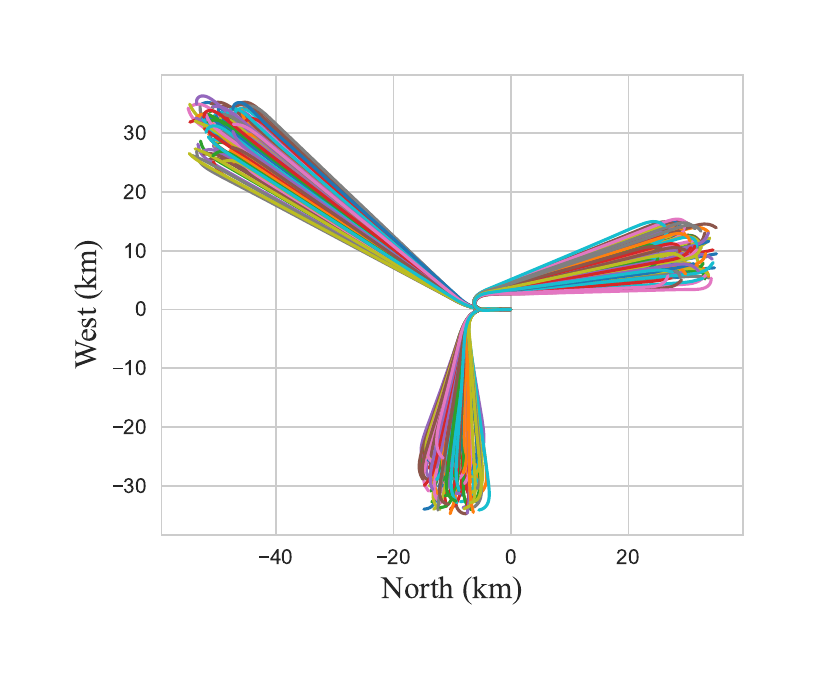}
         \end{subfigure}
         \begin{subfigure}[b]{0.5\textwidth}
         \centering
         \includegraphics[width=0.8\textwidth,trim={2cm 1cm 0cm 2cm},clip]{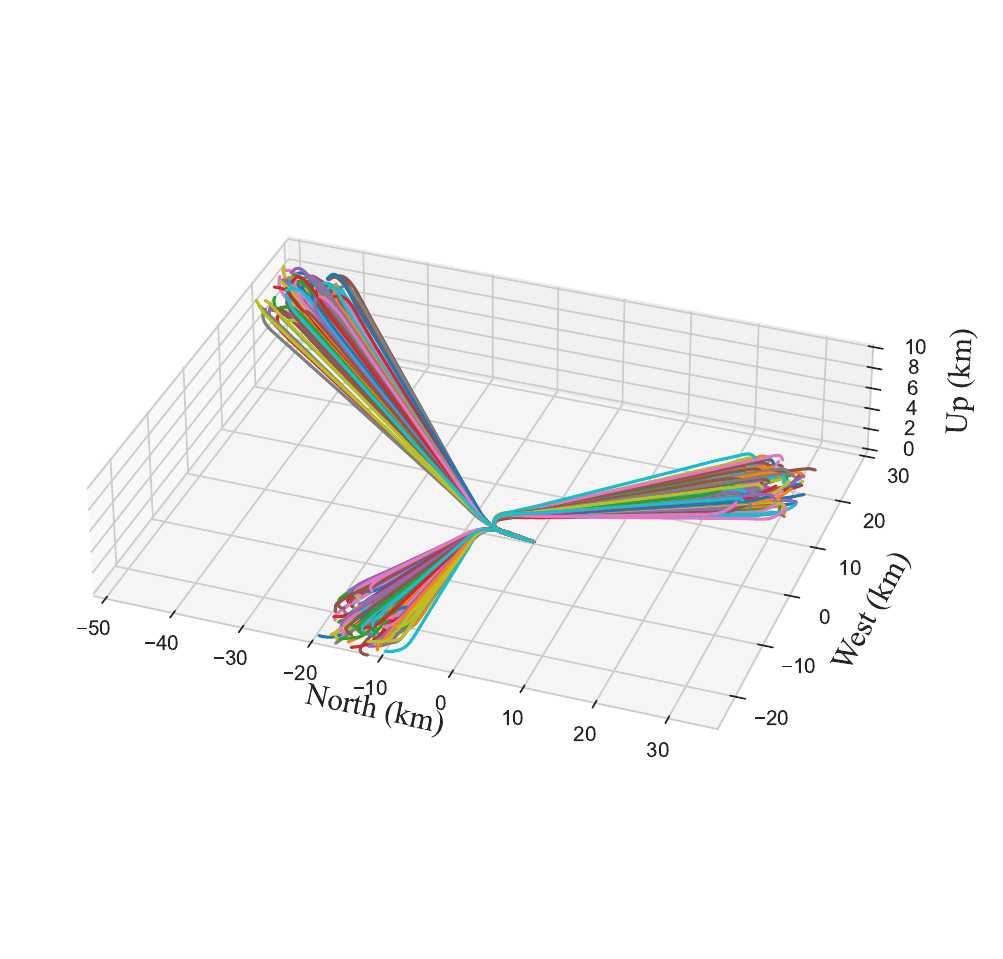}
         \end{subfigure}
\caption{The converged trajectory results with $\gamma = 1$ starting at the initial conditions used for the Monte Carlo study.}
\label{fig:05_traj_result}
\end{figure}

For each initial condition, we categorize the outcome as either a "Success", "Fail by max iteration", or "Fail by divergence". If the first attempt for a given initial condition succeeds, we label it a "Success" and proceed to the next condition. Should the first attempt reach the maximum number of iterations $I_{\max}$, we evaluate the termination criteria $\{(J_{vc} < \epsilon_{vc}) \text{ and } (J_{tr} < \epsilon_{tr})\}$. If any condition in the termination criterion is not satisfied, we double the corresponding weight parameters ($w_{vc}$ and/or $w_{tr}$) and run the second attempt. If the subproblem is not solvable during the iterations of the first attempt, we double $w_{tr}$ and run the second attempt. For the second attempt, a successful convergence results in the "Success" label. However, if the second attempt reaches the maximum number of iterations $I_{\max}$, it is classified as "Fail by max iteration", and if the subproblem is unsolvable in the iterations of the second attempt, we label it as "Fail by divergence". To visualize the chosen initial conditions, we depict the corresponding converged trajectory results with $\gamma = 1$ in Figure \ref{fig:05_traj_result}. For every initial condition and extrapolation parameter, we start with same weight parameters that are $w_{vc}=10^2$ and $w_{tr}=1$ and the maximum number of iterations $I_{\max} = 100$. This outcome categorization is summarized in Algorithm~\ref{alg:outcome_label}.

The simulation results are summarized in Table.~\ref{tab:monte_carlo}. This table provides counts for each label, and for the "Success" label cases, it lists the average cost, average iteration count, and average computation time. Compared to the baseline case with $\gamma = 1$, the algorithm with $\gamma = 1.2$ results in one (0.34\%) additional failed case, but it reduces the average iteration count by about 5 (15.81\%). As the extrapolation is increased further, both the mean cost and the mean iteration count tend to decrease, while the number of failed cases rises. It is worth noting that the mean computation time is generally proportional to the mean iteration count since the extrapolation parameter does not affect the computational complexity of each iteration of the algorithm.

\change{Next, we re-examine the failed cases by repeating the previously described procedure, but this time, starting with weight parameters $w_{vc}=10^3$ and $w_{tr}=20$. For $\gamma = 1$, all 9 out of the 9 failed cases are resolved; for $\gamma = 1.2$, 10 out of the 10 failed cases are resolved; for $\gamma = 1.4$, 10 out of the 11 failed cases are resolved; and for $\gamma = 1.6$, 36 out of the 46 failed cases are corrected. As a result, the proposed algorithm achieves success rates of 100\% for $\gamma=1.0$, 100\% for $\gamma = 1.2$, 99.67\% for $\gamma=1.4$, and 96.67\% for $\gamma=1.6$. In conclusion, the results of the Monte Carlo study show that the proposed method converges reliably for a wide range of initial conditions, and that the xPTR algorithm can accelerate the PTR algorithm, especially with extrapolation parameter values between 1.0 and 1.4.}

\begin{table}[H]
\caption{Monte Carlo simulation results.}
\label{tab:monte_carlo}
\begin{centering}
\begin{tabular}{ccccc}
\hline
$\gamma$ & 1.0 & 1.2 & 1.4 & 1.6\tabularnewline
\hline
Success & 291 (97.0\%) & 290 (96.67\%) & 289 (96.33\%) & 254 (84.67\%) \tabularnewline
\hline
\makecell{Fail by \\ divergence} & 0 (0.0\%) & 0 (0.0\%) & 2 (0.67\%) & 38 (12.67\%) \tabularnewline
\hline
\makecell{Fail by \\  max iteration} & 9 (3.0\%) & 10 (3.33\%) & 9 (3.0\%) & 8 (2.67\%) \tabularnewline
\hline
Mean Cost & 9.9313 & 9.9358 & 9.9182 & 9.6705 \tabularnewline
\hline
\makecell{Mean \\ iteration count} & 28.24 & 23.59 & 21.18 & 21.94 \tabularnewline
\hline
\makecell{Mean \\ computation time (s)} & 13.04 & 10.85 & 10.02 & 11.19 \tabularnewline
\hline
\end{tabular}
\par\end{centering}
\end{table}

\section{Conclusions} \label{sec:conclusion}
We have presented an optimization method for generating aircraft approach and landing trajectories. The problem considers real-world operational constraints including 6-DoF dynamics, a runway alignment constraint, and an obstacle avoidance constraint. To solve the formulated optimal control problem, we develop a novel sequential convex programming algorithm called the extrapolated penalized trust region (xPTR) method, that builds upon the penalized trust region (PTR) algorithm to accelerate convergence by incorporating extrapolation. Through numerical simulations, we show that the proposed method can generate operational trajectories for a wide range of initial conditions.

Planned future work involves considering robust controller synthesis for tracking the generated trajectory under uncertainty \cite{10167750}. One can also consider using the framework proposed in this work in a receding horizon fashion for model predictive control (MPC) \cite{eren2017model}. Closing the loop on the proposed trajectory generator with appropriately designed feedback control laws could also be of interest \cite{garone2017reference}.

\section*{Acknowledgments}
This work was supported by Boeing under Grant 2021-PD-PA-471. We thank the members of the Autonomous Controls Laboratory (ACL) at the University of Washington, especially Purnanand Elango for the discussion and valuable feedback. We also thank Dragos Margineantu for valuable guidance and helpful suggestions.

\bibliography{main}

\end{document}